\documentclass{amsart}
\usepackage{amscd, amssymb}

\def\Kweb#1{
http:\linebreak[2]//www.\linebreak[2]math.\linebreak[2]uiuc.\linebreak[2]edu/\linebreak[2]
{K-theory/\linebreak[2]#1}}
\def\cA{\mathcal A}

\def\cE{\mathcal E}

\def\cJ{\mathcal J}

\def\cO{\mathcal O}

\newcommand{\comment}[1]{}

\newcommand{\Gm}{\mathbb G_m}
\newcommand{\bbH}{\mathbb H}

\newcommand{\N}{\mathbb{N}}

\newcommand{\Q}{\mathbb{Q}}
\newcommand{\R}{\mathbb{R}}
\newcommand{\Z}{\mathbb{Z}}
\newcommand{\frakt}{\mathfrak{t}}

\def\HHF{HH}
\def\inf{\mathrm{inf}}
\def\rht{\mathrm{rht}}
\def\zar{\mathrm{zar}}
\def\harrow{\hookrightarrow}

\newcommand{\bbHz}{\bbH_\zar}

\def\bu{\bullet}

\def\coker{\operatorname{coker}}

\def\hom{\operatorname{hom_\ast}}
\def\Hom{\operatorname{Hom}}

\def\Prim{\operatorname{Prim}}
\def\pro{\operatorname{pro--\!\relax}}

\def\Spec{\operatorname{Spec}}

\def\Tot{\operatorname{Tot}}
\def\map#1{{\buildrel #1 \over \lra}}
\def\smap#1{\,{\buildrel #1 \over \to}\, }
\def\lra{\longrightarrow}
\DeclareMathOperator*{\holim}{holim}
\DeclareMathOperator*{\hocolim}{hocolim}

\newcommand{\SchF}{\mathrm{Sch}/F}

\input xypic
\xyoption{all} 
\numberwithin{equation}{section}

\theoremstyle{plain}
\newtheorem{thm}[equation]{Theorem}
\newtheorem{cor}[equation]{Corollary}
\newtheorem{lem}[equation]{Lemma}
\newtheorem{prop}[equation]{Proposition}

\newtheorem{Cthm}[equation]{Cathelineau's Theorem} 

\theoremstyle{definition}
\newtheorem{defn}[equation]{Definition}
\theoremstyle{remark}
\newtheorem{rem}[equation]{Remark}
\newtheorem{ex}[equation]{Example}
\newtheorem{subremark}{Remark} [equation] 

\begin{document}
\bibliographystyle{plain}

\title[Infinitesimal cohomology and the Chern character]
{Infinitesimal cohomology and the Chern character to negative cyclic homology}

\author{G. Corti\~nas}
\thanks{Corti\~nas' research was partially supported by FSE and by grants
ANPCyT PICT 03-12330, UBACyT-X294, JCyL VA091A05, and MEC MTM00958.
}
\address{Dep. Matem\'atica, FCEyN-UBA \\ Ciudad Universitaria Pab 1\\
1428 Buenos Aires, Argentina\\ and  Dep. \'Algebra\\ Fac. de Ciencias\\
Prado de la Magdalena s/n\\ 47005 Valladolid, Spain.}
\email{gcorti@agt.uva.es}\urladdr{http://mate.dm.uba.ar/\~{}gcorti}

\author{C. Haesemeyer}
\address{Dept.\ of Mathematics, University of Illinois, Urbana, IL
61801, USA} \email{chh@math.uiuc.edu}

\author{C. Weibel}
\thanks{\hskip-8pt Weibel's research was 
supported by NSA grant MSPF-04G-184, and 
by the Oswald Veblen Fund.}
\address{Dept.\ of Mathematics, Rutgers University, New Brunswick,
NJ 08901, USA} \email{weibel@math.rutgers.edu}

\date{\today}

\begin{abstract} There is a Chern character from $K$-theory to
negative cyclic homology. We show that it preserves the decomposition
coming from Adams operations, at least in characteristic zero.
\end{abstract}

\maketitle
\vspace{-25pt} 
\section*{Introduction}

Shortly after the discovery of cyclic homology, Loday and others
raised the question as to whether the Chern character $ch:K_n(A)\to
HN_n(A)$ is compatible with the Adams operations $\psi^k$;
see \cite[0.4.1]{GW94}. Indeed, the respective Adams operations are
defined in very different ways: for $K$-theory in \cite{Q-ICM}
\cite{Kratzer} (see \cite{Soule}), for cyclic homology in \cite{FT, L89}
and for negative cyclic homology $HN_*$ in \cite{LodayHC92, WeibelHA94}.

The following theorem, which we also prove for schemes of finite type over
a field in Theorem \ref{thm:mth} below, answers this in the affirmative for
commutative rings $A$ containing $\Q$.

\begin{thm}\label{main-intro}
For any commutative $\Q$-algebra $A$, the map $ch:K_n(A)\to HN_n(A)$
satisfies $\psi^k ch(x)=k\, ch(\psi^kx)$.
That is, $ch$ sends $K_n^{(i)}(A)$ to $HN_n^{(i)}(A)$.
\end{thm}

Because $K^{(i)}(A)$ and $HN^{(i-1)}(A)$ are the eigenspaces for
$\psi^k=k^i$ (see \cite[4.5]{LodayHC92}), the two formulations in this
theorem are equivalent. The shift in indexing for $HN$ arises from
the desire to have $\lambda^k=(-1)^{k-1}k^i$ on $HN^{(i)}(A)$;
see \cite[4.5.4]{LodayHC92}.

\begin{ex}\label{ex:ch(t)}
To see the difference in eigenvalues, consider the ring
$R = F[t,t^{-1}]$. It is well known that the element $t\in K_1(R)$
satisfies $\psi^k(t)=k\cdot t$ for all $k$, while its character in
$HN_1(R)$ satisfies $\psi^k ch(t)=k^{2}ch(t)$; see \cite[8.4.7]{LodayHC92}.
The image of $ch(t)$ in $HH_1(R)=\Omega^1_R$ is $dt/t$, and again
$\psi^k(dt/t)=k^2(dt/t)$.
\end{ex}

Several special cases of this theorem have been addressed in the literature.
The case of nilpotent ideals was asserted in the proof of \cite[7.5.5]{FT}
and established by Cathelineau in \cite{C91} (see our Appendix \ref{app:A}).
The case $K_0(A)\to HN_0(A)$, as well as the case $n<0$, was settled in
\cite{Weibel93}, since $HN_n(A)\to HP_n(A)$ is an isomorphism for $n\le0$
and the composition $K_n(A)\to HN_n(A)\to HP_n(A)$ was shown to have this
property in {\it loc.\ cit.}
The Dennis trace map, which is the composition
$K_n(A)\to HN_n(A)\to HH_n(A)$, was shown to be compatible
with the Adams operations by Kantorowitz in \cite{Kwz99}.
Given the results of \cite{chsw}, Gillet and Soul\'e proved a very
similar result in \cite[3.2.2 and 6.1]{GS99}, using
the universal total Chern class.

Our method is to use infinitesimal cohomology to reduce the problem to
the nilpotent case considered by Cathelineau. To that end,
the goal of the first three sections is to interpret the infinitesimal
cohomology $H^*(X_\inf,HN)$ in terms of $H^*(X_\inf,\cO)$.
In section \ref{sec:infh}, we recall some basic facts about
infinitesimal cohomology, introduced by Grothendieck in \cite{Dix}.
Section \ref{sec:pro-homology} reviews some elementary facts about
pro-homological algebra which we need in section \ref{williefix}
to prove a pro-version of the Hochschild-Kostant-Rosenberg Theorem.
The interpretation of $H^n(X_\inf,HN)$, and the fact that it vanishes
for $n<0$, occurs in \ref{ex:hinfHN}.

Section \ref{sec:Hinf} introduces sheaf hypercohomology spectra for
the infinitesimal topology, as a generalization of Grothendieck's
construction. These ideas are applied to the $K$-theory spectrum in
section \ref{sec:Kinf}, where the space $\bbH(X_\inf,K)$ is compared to
the fiber $K^{\inf}$ of the Chern character. Here we establish the
technical fact that $HN_n(X)\to\pi_{n-1}\bbH(X_\inf, HNI)$ is an isomorphism
for $n\le-2$, and is injective for $n=-1$, where the presheaf $HNI$
on $X_\inf$ is defined in \ref{def:KIT} so that $\bbH(X_\inf, HNI)$
is the fiber of $\bbH(X_\inf,HN)\to HN(X)$.

Section \ref{sec:Cath-schemes} extends Cathelineau's result to schemes,
and to infinitesimal hypercohomology.
Finally, our main theorem is proven in section \ref{sec:main} as a
special case of Theorem \ref{thm:mth}, which in turn follows from the
technical fact mentioned above.

In Appendix \ref{app:A}, we provide a technical correction to the preprint
of \cite{OW}, and to its use in the
proof of Cathelineau's theorem in \cite{C91}.  Although this
correction is well known to the experts, it has not appeared in print
before.  In Appendix \ref{app:B}, we give a simplicial presheaf version
(and a spectrum version) of the same theorem. It is this version of
Cathelineau's construction that we need in order to prove our main theorem.
Appendix \ref{app:B} depends upon some technical results about the model
structure of simplicial presheaves (of sets); these results
are proven in Appendix \ref{app:C}.

\medskip

\noindent{\bf Notation.} 
We shall write $\SchF$ for the category of schemes essentially of finite
type over a field $F$. Objects of $\SchF$ shall be called {\it $F$-schemes}.
If $k\subset F$ is a subfield, we write $\Omega^p_{/k}$ for the sheaf of
$p$-differential forms; we will write $\Omega^p$ when $k$ is clear
from the context.
If $H$ is a functor on $\SchF$ and $X=\Spec(A)$, we shall sometimes write
$H(A)$ instead of $H(\Spec A)$; for example, $H^*(A,\Omega^p)$ is used for
$H^*(\Spec A,\Omega^p)$.

We use cohomological indexing for all chain complexes in this paper;
for a complex $C$, $C[p]^q = C^{p+q}$.
For example, the Hochschild, cyclic, periodic and negative cyclic
homology of schemes over a field $k$ (such as $F$-schemes over a field
$F\supseteq k$) can be defined using the Zariski hypercohomology of
certain presheaves of complexes; see \cite{WeibelHC} and
\cite[2.7]{chsw} for precise definitions.  We shall write these
presheaves as $HH(/k)$, $HC(/k)$, $HP(/k)$ and $HN(/k)$,
respectively, omitting $k$ from the notation if it is clear from the context.

If $\cE$ is a presheaf of spectra on $\SchF$ (or just on $X$), we
write $U\mapsto \bbHz(U,\cE)$ for Thomason's sheaf hypercohomology
spectrum \cite[1.33]{AKTEC}. Jardine showed in \cite[3.3]{JardineSPS}
that $\bbHz(-,\cE)$ is the fibrant replacement for $\cE$ in the model
structure of \cite[1.5]{JardineSHT}; see \cite[D.5]{TT}.
We say that $\cE$ satisfies {\it Zariski descent} on $\SchF$ (or on
$X$) if the natural maps $\cE(U)\to \bbHz(U,\cE)$ are homotopy
equivalences for all $U$ in $\SchF$ (resp., $U\subset X$).

\smallskip
It is well known (see \cite[10.9.19]{WeibelHA94})
that there is an Eilenberg-Mac\,Lane functor $C\mapsto|C|$
from chain complexes of abelian groups to spectra, and from presheaves
of chain complexes of abelian groups to presheaves of spectra.
This functor sends quasi-isomorphisms of complexes to weak
homotopy equivalences of spectra, and satisfies $\pi_n(|C|)=H^{-n}(C)$;
the loop space is $\Omega|C|\simeq |C[-1]|$.
For example, applying $\pi_n$ to the Chern character $K\to|HN|$ yields maps
$K_n(A)\to H^{-n}HN(A)=HN_n(A)$.
\newpage
\section{Infinitesimal cohomology of sheaves}\label{sec:infh}

Recall from \cite{Dix} that a closed immersion of schemes $U\harrow T$
is called a {\it thickening} of $U$ if its ideal of definition is
nilpotent; by abuse of notation, we write $T$ to mean $U\harrow T$.  For
$X\in\SchF$, the {\it infinitesimal site} $X_{\inf}$ consists of the
category $\inf(X)$ and its coverings, which we now define. Objects of
$\inf(X)$ are thickenings $U\harrow T$, where $U$ is an open subscheme
of $X$; morphisms from $U\harrow T$ to $U'\harrow T'$ are morphisms
$T\to T'$ in $\SchF$ under inclusions $U\subseteq U'$. A covering of $T$ is a
family of morphisms $\{T_i\to T\}$ such that the $T_i$ form a Zariski
open covering of $T$.

A sheaf $\cE$ on $X_{\inf}$ is the same thing as a compatible
collection of Zariski sheaves $\{\cE_T\in Sh(T):T\in \inf(X)\}$
(see \cite[\S5]{BO}, \cite[4.1]{Dix}). It follows that $X_\inf$ has
enough points, namely the Zariski points on the thickenings $T$.

Note that by definition, the forgetful functor
\[
u: \inf(X)\to \SchF,\ \ (U\harrow T)\mapsto T,
\]
is a morphism of sites $(\SchF)_\zar \to X_\inf$.
Thus if $\cE$ is a Zariski sheaf on $\SchF$, its restriction to
$\inf(X)$ defines a sheaf $u^*\cE$ on $X_\inf$;
by abuse, we will write $\cE$ for $u^*\cE$.
The usual global sections functor takes a sheaf $\cE$ on $X_\inf$ to:
\[
H^0(X_\inf,\cE) = \Gamma_\inf(\cE) = \lim_{T\in\inf(X)}H^0(T,\cE_T),
\]
and the infinitesimal cohomology of a sheaf of abelian groups is
defined as:
\begin{equation}\label{def:Hinf1}
H^*(X_\inf,\cE)=H^*\ \bbH(X_\inf,\cE), \quad
\text{where }\bbH(X_\inf,\cE)=\R\Gamma_\inf(\cE).
\end{equation}

\begin{ex}\label{ex:338}
For affine $X=\Spec(A)$, we can compute infinitesimal cohomology using
the method outlined by Grothendieck in \cite[p.~338]{Dix}.  Suppose
for simplicity that $\cE$ is a quasi-coherent sheaf on $\SchF$, and
that the natural map $j^*\cE_{T'}\to \cE_T$ is an isomorphism for
every closed embedding $j:T\hookrightarrow T'$ in $\inf(X)$; in this
case $\cE$ is called a {\it crystal}; see \cite[2.12]{BO}.

Let $A=S/I_1$ be a presentation of $A$ as a quotient of a smooth
$F$-algebra $S$, and write $I_\nu$ for the kernel of
$S^{\otimes\nu}\to A$.  Then each $S^{\otimes\nu}_\bullet=
\{ S^{\otimes \nu}/I_\nu^m\}_m$ is a tower of algebras, and
$[\nu]\mapsto \{ S^{\otimes\nu+1}_\bullet\}$ is a
cosimplicial tower of algebras. Following \cite{Dix}, set
$Y^\nu_m=\Spec (S^{\otimes \nu}/I_\nu^m)$.
Further, define $\cE(Y^{\nu}_\bullet)= \varprojlim_m \cE(Y^\nu_m)$.
Regarding the cosimplicial group
$\cE(Y^{\ast+1}_\bu)$ as a cochain complex, we have:
\[ H^*(X_\inf,\cE) = H^* \cE(Y^{\ast+1}_\bullet).  \]
\end{ex}

Here is a modern interpretation of Grothendieck's argument in
\cite[5.2]{Dix}.  The functor $\Gamma_\inf$ factors as the forgetful
functor $\Gamma_\zar$ from sheaves to presheaves, followed by the
inverse limit functor. The forgetful functor preserves injectives as
it is right adjoint to sheafification. Following \cite[XI.6]{BK} and
\cite[5.32]{AKTEC}, we write $\holim_T$ for $\R\lim_{T\in\inf(X)}$ and
$\bbHz(T,\cE)$ for $\R\Gamma_\zar(T,\cE)$.
It follows from \cite[10.8.3]{WeibelHA94} that we have an isomorphism
in the derived category:
\begin{equation}\label{sheafHinf}
\bbH(X_\inf,\cE) \simeq \holim\nolimits_T \bbHz(T,\cE).
\end{equation}

\begin{lem}\label{lem:5.1}
Suppose that $X=\Spec(A)$ and $A=S/I$ 
for a smooth algebra $S$.
Then the simplicial cotower
$[\nu]\mapsto \{ Y^{\nu+1}_\bullet\}$ 
of \ref{ex:338} is right cofinal in $\inf(X)$.
\end{lem}

Right cofinality means that for each object $T$ of $\inf(X)$,
the category $T/i$ is contractible, where $i$ denotes the inclusion of
the simplicial cotower into $\inf(X)$.

\goodbreak
\begin{proof}

By  \cite[5.1]{CortInf}, the simplicial cotower is right cofinal in
the subcategory $X/\inf(X)$ of thickenings of $X$.
We assert that the proof of \cite[5.1]{CortInf} goes through
{\it mutatis mutandis}.
Given a thickening $U\hookrightarrow T$, set $B=\cO(T)$, $C=\cO(U)$ and pick
$r$ so that $\ker(B\to C)^r=0$. Consider the cosimplicial object
$h_r:\Delta\to\inf(X)^{op}$ defined by the $Y^{\nu+1}_r$.

As in {\it op.\ cit.}, the inclusion $T/h_r\subset T/i$ is a homotopy
equivalence, $(T/h_r)^{op}$ is the homotopy colimit of the discrete
simplicial set $\Hom_{\inf(X)}(T,Y^{\ast+1}_r)$, and this simplicial set
is contractible. Since $T/h_r$ is contractible, so is $T/i$.
\end{proof}

Let $[\nu]\mapsto C_\nu$ be a cosimplicial abelian group or,
more generally, a cosimplicial complex of abelian groups.
Recall from \cite[5.32]{AKTEC} that the Bousfield-Kan total complex
$\holim_{\nu\in\Delta}C^*_\nu =
\Tot_\nu C^*_\nu$ is just a specific Cartan-Eilenberg resolution
of the associated total cochain complex $C^*$.

\begin{prop}\label{sheaf-holim}
Suppose that $X=\Spec(A)$ and that $\cE$ is a cochain complex of
quasi-coherent sheaves on $\SchF$.
\begin{enumerate}
\item $\holim_T \cE(T) \simeq \Tot_\nu(\R\varprojlim_m)~\cE(Y^{\nu+1}_m)$.

\item  Assume that $j^*\cE_{T'}\to \cE_T$ is onto for
every closed embedding $j:T\hookrightarrow T'$ in $\inf(X)$.
Writing $\cE(Y^{\nu}_\bu)$ for $\varprojlim_m \cE(Y^\nu_m)$
as in Example \ref{ex:338}, we have:
$\bbH(X_\inf,\cE) \simeq \Tot_\nu \cE(Y^{\nu+1}_\bu)$.
\end{enumerate}
\end{prop}
\goodbreak

\begin{proof} By Lemma \ref{lem:5.1} and the
cofinality theorem \cite[XI.9.2]{BK}, $\holim_T$ is equivalent to
$\holim_{\N\times\Delta}=\holim_{\Delta}\holim_{\N}=\Tot(\R\varprojlim_m)$.
Part (1) follows because, by hypothesis, $\cE(T)=\bbHz(T,\cE)$ for
all affine $T$. Under the assumptions of (2),
$\R\varprojlim_m \cE(Y^\nu_m)$ can be
replaced by $\varprojlim_m \cE(Y^\nu_m)=\cE(Y^\nu_\bu)$.
\end{proof}

\begin{subremark}\label{rem:notcrystal}
Let $k\subset F$ be a subfield. The sheaf $\Omega^p=\Omega^p_{/k}$
of $p$-differential forms is not a crystal for $p\ne0$ in the sense of
\ref{ex:338}. However, it satisfies the hypotheses of \ref{sheaf-holim}(2),
allowing us to compute $\bbH(X_\inf,\Omega^p)$ using $\Tot_\nu$.
\end{subremark}

\smallskip
Let $\cE$ be a cochain complex of sheaves of abelian groups on $\inf(X)$.
\begin{lem}\label{infZdescent}
The presheaf $V\mapsto \bbH(V_\inf,\cE)$ satisfies
Zariski descent on $X$.
\end{lem}

\begin{proof}
For every Zariski open $j:V \subset X$, the inclusion
$\inf(V)\subset\inf(X)$ has a
right adjoint $\rho:\inf(X)\to\inf(V)$, sending thickenings of $U$ to
thickenings of $U\cap V$. For simplicity, we shall write $T\cap V$
for the thickening of $U\cap V$ corresponding to $T$, so that
$\rho(U\hookrightarrow T)$ is $U\cap V\hookrightarrow T\cap V$.
By \ref{lem:5.1},
the Bousfield-Kan cofinality theorem \cite[XI.9.2]{BK}
applies to say that the map of Cartan-Eilenberg resolutions
\addtocounter{equation}{-1}
\begin{subequations}
\begin{equation}\label{eq:j*E}
\holim_{T\in\inf(X)}\bbHz(T\cap V,\cE) \to \holim_{\inf(V)}\bbHz(T,\cE)
\end{equation}
\end{subequations}
is a quasi-isomorphism.
Now if $\{V_1,V_2\}$ is a Zariski cover of an open $V$, we have a
distinguished triangle 
for every $T$ in $\inf(V)$:
\[
\bbHz(T,\cE) \to \bbHz(T\cap V_1,\cE)\times\bbHz(T\cap V_2,\cE) \to
\bbHz(T\cap V_1\cap V_2,\cE).
\]
Taking the homotopy limit over $T$, we obtain the triangle 
\[
\bbH(V_\inf,\cE) \to \bbH(V_1{}_\inf,\cE)\times\bbH(V_2{}_\inf,\cE) \to
\bbH((V_1\cap V_2)_\inf,\cE),
\]
which implies the assertion that $\bbH(-_\inf,\cE)$
satisfies Zariski descent on $X$.
\end{proof}

\begin{subremark}
Formula \eqref{eq:j*E} says that for any Zariski open
$j:V\subset X$ and any complex of presheaves of abelian groups
$\cE$ on $\inf(V)$, the direct image
$j_*\cE(T)=\cE(T\cap V)$ satisfies
$\bbH(V_\inf,\cE)\cong\bbH(X_\inf,j_*\cE)$.
\end{subremark}
\goodbreak
Here is an extension of Grothendieck's theorem \cite[4.1]{Dix}:

\begin{thm}\label{thm:infOmega}
The brutal truncations $\Omega^*\to\Omega^{\le i}\to\Omega^0=\cO$
induce homotopy equivalences
\[
\bbH(X_\inf,\Omega^*)\simeq \bbH(X_\inf,\Omega^{\le i})
\simeq \bbH(X_\inf,\cO).
\]
\end{thm}

\begin{proof}
From the hypercohomology spectral sequence (see \cite[5.7.9]{WeibelHA94})
$$E_1^{p,q}=H^q(X_\inf,\Omega^p)\Rightarrow H^{p+q}(X_\inf,\Omega^*),$$
we see that it suffices to prove that $\bbH(X_\inf,\Omega^p)=0$ for
$p>0$. By induction on the size of a separated cover, using Lemma
\ref{infZdescent}, we are reduced to the case in which $X$ is separated.
A similar induction reduces us to the case in which $X$ is affine.

The affine  case was established in \cite[7.9]{CortDR} using the method of
Example \ref{ex:338} above. 
As argued in {\it loc.\ cit.}, it suffices to show that
$\{\Omega^p_{S^{\otimes\bullet}/I_\bullet}\}_\nu$ is acyclic, where $S$
is a symmetric algebra. This complex is the inverse limit of the tower
$$
\{\Omega^p_{S^{\otimes\bullet}}/I_\bullet^m\Omega^p_{S^{\otimes\bullet}}\}_\nu.
$$
When $I=0$, a more or less explicit chain contraction is given in
{\it loc.\ cit.}, and is defined by differential operators.

Since a contraction of this kind must be continuous for the
$I_\bullet$-adic topology, the result follows.
\end{proof}

Associated to the de\,Rham mixed complex $(\Omega_R,0,d)$ is the
Connes double complex $HC\Omega_R$ for cyclic homology; see
\cite[2.5.10]{LodayHC92}, \cite[9.8.8]{WeibelHA94}. We can also form
double complexes $HP\Omega_R$ and $HN\Omega_R$ for periodic and negative
cyclic homology. The following is immediate from Theorem \ref{thm:infOmega}
and the hypercohomology spectral sequence for the row filtration
on these double complexes.

\begin{cor}\label{cor:HComega}
For any $X$ in $\SchF$,
$\bbH(X_\inf,HC\Omega)\simeq \prod_{i\ge0}\bbH(X_\inf,\cO)[2i]$.

Similarly, for the periodic and negative cyclic variants we have
\[
\bbH(X_\inf,HP\Omega)\simeq \prod_{i\in\Z}\bbH(X_\inf,\cO)[2i]
\quad\text{and}\quad
\bbH(X_\inf,HN\Omega)\simeq \prod_{i\le0}\bbH(X_\inf,\cO)[2i].
\]
\end{cor}

\bigskip
\section{pro-homological algebra}\label{sec:pro-homology}

In order to prove the main result (Theorem \ref{thm:pro-hkr})
in the next section, we need some
elementary results on the homological algebra of pro-objects.
No great originality is claimed for the results in this section.

\smallskip
We recall from \cite[A.4.5]{AM} that if $\cA$ is an abelian category
then the pro-category $\pro\cA$ is also abelian, and so is the full
subcategory of $\pro\cA$ consisting of towers $\{A_m\}=\{A_m\}_{m}$
(indexed by the natural numbers $m\ge0$).

A tower $\{A_m\}$ is isomorphic to $0$
in $\pro\cA$ if and only if it satisfies the trivial Mittag-Leffler
condition that for every $m$ there exists a $j>m$ such that $A_j\to A_m$
is zero. A strict morphism of towers $\{f_m\}:\{A_m\}\to\{B_m\}$ is
an isomorphism in $\pro\cA$ if and only if the kernel $\{\ker(f_m)\}$
and cokernel $\{\coker(f_m)\}$ are isomorphic to $0$ in $\pro\cA$.

If $\{i_m\}$ is a strictly increasing sequence of natural numbers
then $\{A_{i_m}\}\to\{A_m\}$ is an isomorphism in $\pro\cA$, and every
morphism $\{A_i\}\to\{B_m\}$ in $\pro\cA$ is represented by a strict
morphism of towers $\{A_{i_m}\}\to\{B_m\}$.
\goodbreak

If $R$ is a ring, $M$ an $R$-module, and $J\subset R$ an ideal, we
write $M/J^\infty M$ for the pro-$R$-module $\{M/J^mM\}_m$.

\begin{lem}\label{lem:pro-exact}
Let $R$ be a noetherian ring and $J$ an ideal. Then
$M\mapsto M/J^\infty M$ is an exact functor from the category of
finitely generated $R$-modules to pro-$R$-modules.
\end{lem}

\begin{proof} Suppose that $0\to L\to M\to N\to0$ is an exact sequence
of finitely generated $R$-modules. Then $M/J^\infty M\to N/J^\infty N$
is onto with kernel  $\{L/J^mM\cap L\}$. By the Artin-Rees lemma,
$\{J^mL\}\to \{J^mM\cap L\}$ is a pro-isomorphism. It follows that
$\{L/J^mL\}\to\{L/J^mM\cap L\}$ is also a pro-isomorphism.
\end{proof}

\begin{ex}\label{ex:pro-Omega}
If $I$ is an ideal in a ring $S$, then the pro-$R$-modules
$\Omega^p_S/I^\infty\Omega^p_S$ and $\{\Omega^p_{S/I^m}\}$ are isomorphic
for all $p$, with the isomorphism coming from the Fundamental Exact Sequence

$I/I^2\to\Omega_S/I\Omega_S\to\Omega_{S/I}\to0$. In particular,
$\varprojlim\Omega^p_S/I^m\Omega^p_S=\varprojlim\Omega^p_{S/I^m}$.
\end{ex}

\begin{lem}\label{lem:pro-pro}
Let $R$ be a ring, $J\subset R$ an ideal and $f:M\to N$ a surjective
homomorphism of $R$-modules. Let $P:=\{P_m\}$ be a pro-$R$-module
such that each $P_m$ is a projective $R/J^m$-module.
Given any morphism $g:P\to N/J^\infty N$ of pro-$R$-modules,
there exists a pro-homomorphism $h:P\to M/J^\infty M$ making the
following diagram commute.
\[
\xymatrix{&P\ar@{.>}[dl]_h\ar[d]^g\\
          M/J^\infty M\ar[r]_f& N/J^\infty N}
\]
\end{lem}

\begin{proof}
Choose a representative $\{g_m:P_{i_m}\to N/J^mN\}$ of $g$.
Then each $Q_m=P_{i_m}/J^mP_{i_m}$ is a projective $R/J^m$-module,
and there is a pro-module isomorphism $\{Q_m\}\cong\{P_m\}$.

Replacing $\{P_m\}$ by $\{Q_m\}$ if
necessary, we may assume that $\{g_m:P_m\to N/J^mN\}$ is a strict map.

We will construct a strict lift $\{h_m\}$ by induction on $m$. The case
$m=1$ is clear because $P_1$ is a projective $R/J$-module.
Inductively, we have a lift $h_m:P_m\to M/J^mM$ of $g_m$ and hence
a map $P_{m+1}\to M/J^mM\times_{N/J^mN}N/J^{m+1}N$. Since
$M/J^{m+1}\to M/J^mM\times_{N/J^mN}N/J^{m+1}N$ is onto and $P_{m+1}$
is projective, we get the desired lift $P_{m+1}\to M/J^{m+1}$
compatible with $h_m$ and $g_{m+1}$.
\end{proof}

\begin{cor}\label{cor:pro-acyc}
Let $R$ be a noetherian ring,  $M$ a finitely generated $R$-module, and
$\epsilon:L_*\to M$ a resolution by finitely generated $R$-modules.
Let $J\subset R$ be an ideal. Then $L_*/J^\infty L_* \to M/J^\infty M$ is
a resolution in the category of pro-$R$-modules.

If in addition $P_*=\{P_{*,m}\}_m$ is a 
chain complex of pro-$R$-modules, such that each $P_{n,m}$ is a
projective $R/J^m$-module, and $P_{*}\to M/J^\infty M$ is a chain map of
pro-modules, then there is a chain map of pro-modules $P_*\to L_*/J^\infty L$,
unique up to pro-chain homotopy, making the following diagram commute.
\[
\xymatrix{&P_*\ar@{.>}[dl]\ar[d]\\ L_*/J^\infty L \ar[r]&M/J^\infty M}
\]
\end{cor}

\begin{proof} The first assertion is just Lemma \ref{lem:pro-exact}.
Given this, the proof of the usual Comparison Theorem
\cite[2.2.6]{WeibelHA94} goes through using Lemma \ref{lem:pro-pro}.
\end{proof}

\begin{lem}\label{lem:pro-nul}
Let $R$ be a ring and let $r_1,\dots,r_n\in R$
be such that $R=(r_1,\dots,r_n)R$.
Then for every pro-$R$-module $M=\{M_m\}$, $M$ is pro-isomorphic to $0$
if and only if each $M[1/r_i]$ is.
\end{lem}

\begin{proof} Straightforward, using the trivial Mittag-Leffler condition.
\end{proof}

If $\{P_*=P_{\ast,m}\}$ is a pro-chain complex, meaning that $P_{\ast,m}$
is a chain complex for each $m$ and the $P_{\ast,m+1}\to P_{\ast,m}$ are
chain maps, then we may regard $P$ as a chain complex of pro-objects.
Conversely, any bounded below chain complex of pro-objects is pro-isomorphic
to a pro-chain complex via the re-indexing trick.

\begin{lem}\label{pro-homo}
Let $R$ be a ring, $J\subset R$ an ideal, $P_*$ and $Q_*$ two bounded below
pro-complexes of $R$-modules, and $f:P_*\to Q_*$ a homomorphism of complexes
of pro-$R$-modules. Assume that at each level $m$ and each degree $n$,
both $P_{n,m}$ and $Q_{n,m}$ are $R/J^m$-modules, that
$P_{n,m}$ is projective, and that each complex $Q_{*,m}$ is acyclic.

Then for every additive functor $F$ from $R$-modules to an abelian category
$\cA$, and each $n$, the induced map $H_n(F(P_*))\to H_n(F(Q_*))$

is zero in $\pro\cA$.
\end{lem}

\begin{proof}
We claim that by re-indexing $P_*$ we may assume that
$f$ is a strict chain map of towers. In this case, because each $P_{*,m}$
is a complex of projectives and each $Q_{*,m}$ is acyclic, each $f_m$
is chain-homotopic to zero by the usual Comparison Theorem
\cite[2.2.6]{WeibelHA94}, and the result follows.

To see the claim, note that

we can choose a function $h:\N^2\to\N$,
strictly increasing with respect to each variable separately,

together with a representative $f_{n,m}:P_{n,h(n,m)}\to Q_{n,m}$,
so that the following diagram commutes for all $n$ and $m$:
\[
\xymatrix{P_{n+1,h(n+1,m+1)}\ar[d]_{\partial}\ar[r] &
P_{n+1,h(n+1,m)}\ar[d]_{\partial}\ar[r]^{~f_{n+1,m}}&
Q_{n+1,m}\ar[d]^\partial\\
  P_{n+1,h(n+1,m)}\ar[r] & P_{n,h(n,m)}\ar[r]^{f_{n,m}}& Q_{n,m}.}
\]
(We have abused notation by omitting all notation for transition maps

$P_{*,j}\to P_{*,i}$.)
Set $P'_{n,m}=P_{n,h(n,m)}/J^mP_{n,h(n,m)}$; this is a projective
$R/J^m$-module. Because $P_*$ is a pro-complex, $P'_*$ is a
tower of chain complexes, pro-isomorphic to $P_*$,
and the induced map $P'_*\to Q_*$ is a strict map of towers of complexes.
\end{proof}

\bigskip
\section{Pro-Hochschild-Kostant-Rosenberg theorem}\label{williefix}

Let $F$ be a field, $S$ an algebra over $F$, and $I\subset S$ an ideal.
Write $\Omega$ and $\HHF$ for $\Omega^*_{/F}$ and Hochschild homology
taken over $F$, respectively. Then the shuffle product
\cite[9.4.4]{WeibelHA94} induces a map of graded pro-$S$-modules:

\begin{equation}\label{pro-hkr}
\{\Omega^p_{S/I^m}\}_m\to \{\HHF_p(S/I^m)\}_m.
\end{equation}

The purpose of this section is to prove the following theorem
for fields $F$ of arbitrary characteristic,
and its analogue for subfields $k\subseteq F$ in characteristic~0.

\begin{thm}\label{thm:pro-hkr}
If $S$ is essentially of finite type and smooth over $F$ then
\eqref{pro-hkr}
is a pro-isomorphism for every ideal $I\subset S$.
\end{thm}

The proof of Theorem \ref{thm:pro-hkr} which we shall give is an
adaptation to the pro-setting of the proof of the
Hochschild-Kostant-Rosenberg theorem \cite{hkr} \cite[9.4.7]{WeibelHA94}.
In the proof, all vector spaces, tensor products, algebras, and
differential forms will be taken over $F$.

\addtocounter{equation}{-1}
\begin{subequations}

\begin{proof} 
Let $\cJ=\ker(S^e\to S)$, where $S^e=S\otimes S$.
Because $\cJ/\cJ^2=\Omega^1_{S}$ is projective,
$\cJ$ is locally a complete intersection. Thus we can find
elements $\alpha_1,\dots,\alpha_n\in S^e$
such that $\Spec(S)=V(\cJ)\subset\bigcup_i\Spec S^e[1/\alpha_i]$
and such that for
each $i$, the ideal $\cJ[1/\alpha_i]\subset S^e[1/\alpha_i]$ is
a complete intersection. Let $s_i$ be the image of $\alpha_i$ under
the map $S^e\to S$; we have
$S^e[1/\alpha_i]/\cJ[1/\alpha_i]=S[1/s_i]$, and $\Spec
S=\bigcup_i\Spec S[1/s_i]$. Upon replacing $\alpha_i$ by
$(s_i\otimes s_i)\alpha_i$ if necessary, we may assume that
$S^e[1/\alpha_i]$ is a localization of $S[1/s_i]^e$. Lemma
\ref{lem:pro-nul} applied to the pro-$S$-modules given by the kernel
and cokernel of the map \eqref{pro-hkr} shows that it suffices to check
that the latter becomes an isomorphism after inverting each $s_i$.
Fixing $i$, replacing $S$ by $S[1/s_i]$, and setting $\alpha=\alpha_i$,
we may therefore assume that $\cJ_\alpha=\cJ[1/\alpha]$ is a
complete intersection
in $S^e_\alpha=S^e[1/\alpha]$, and that $S=S[1/\alpha]$. Put
\[
J=(S\otimes I+I\otimes S)[1/\alpha] 
\quad\text{and}\quad J^{(m)}=(S\otimes I^m+I^m\otimes S)[1/\alpha].
\]
Note that $J^{(m)}\subset J^m$ and $J^{2m}\subset J^{(m)}$, so
$J^{(\bu)}$ and $J^\bu$ are equivalent filtrations and the pro-algebra map
$\{(S/I^m)^e=S^e/J^{(m)}\}\to\{S^e/J^\infty\}$ is an isomorphism.
\smallskip

Consider the bar resolution $C^{bar}_*(S)\map{\sim}S$
of $S$ as a $S^e$-module; it has
\[
C^{bar}_n(S) = S^e\otimes S^{\otimes n} = S^{\otimes n+2}\qquad (n\ge 0).
\]
Because $S^e_\alpha$ is a localization of $S^e$,
$S^e_\alpha\otimes_{S^e}C^{bar}_*(S)$ is an $S^e_\alpha$-projective
resolution of $S$. Since $S^e_\alpha/J^{(m)}=(S/I^m)^e_\alpha$,
\[
Q_{*,m} = S^e_\alpha \otimes_{S^e} C^{bar}_*(S/I^m)
  = (S/I^m)^e_\alpha\otimes_{S^e}C^{bar}_*(S/I^m)
\]
is a projective $(S/I^m)^e_\alpha$-module resolution of $S/I^m$.
Similarly, its quotient
\[
P_{\ast,m}=S^e_\alpha/J^m\otimes_{S^e}C^{bar}_*(S/I^m)
\]
is a complex of projective $S^e_\alpha/J^m$-modules. Setting
$S^e_\alpha/J^{(\infty)}=\{S^e_\alpha/J^{(m)}\}$ we have an
isomorphism of pro-complexes
\begin{equation}\label{equivfil}
Q_* = S^e_\alpha/J^{(\infty)}\otimes_{S^e}C^{bar}_*(S/I^\infty S)
\map{\sim} S^e_\alpha/J^\infty\otimes_{S^e}C^{bar}_*(S/I^\infty S) = P_*.
\end{equation}

Choose $x_1,\dots,x_d\in\cJ$ whose images in $J_\alpha$ form a regular
sequence of generators, and write $L_*$ for the Koszul complex
$K(S^e_\alpha;x_1,\dots,x_d)$. Then $L_*$ is an
$S^e_\alpha$-projective resolution of $S$, of the form
$L_*=S^e_\alpha\otimes\wedge^*V$, where $V$ is a $d$-dimensional
$F$-vector space with basis $\{v_1,\dots,v_d\}$, and $L_1\to L_0$
sends $v_i$ to $x_i$.

\goodbreak
By the Comparison Theorem, there is a chain equivalence
$L_*\map{\epsilon} S^e_\alpha\otimes C^{bar}_*(S)$ of projective
resolutions of $S$, unique up to chain homotopy.  We can in fact
choose $\epsilon$ to induce the map \eqref{pro-hkr}. To do so, write

$x_i=\sum_j r_{ij}(s_{ij}\otimes 1-1\otimes s_{ij})$,
with $r_{ij}, s_{ij}\in S$, and set
$\epsilon(v_i) = \sum_j r_{ij}\otimes s_{ij}\otimes 1$.
This defines a map $V\to S^e_\alpha\otimes_{S^e}C^{bar}_1(S)$,
and we extend it to an $S^e_\alpha$-linear homomorphism of chain complexes
using the shuffle product of $C^{bar}_*(S)$.

Composing $\epsilon$ with $S^e_\alpha\otimes_{S^e}C^{bar}_*(S)\to P_{*,m}$
induces maps $L_*/J^mL_* \to P_{*,m}$ and hence a strict map
of pro-complexes
\begin{equation}\label{pro-epsilon}
\epsilon:L/J^\infty L\to P_*
\end{equation}
which covers the identity of $S/I^\infty S$.
Since $S/I^m=S\otimes_{S^e_\alpha} S^e_\alpha/J^m$, tensoring
\eqref{pro-epsilon} 
over $S^e_\alpha$ with $S$ and using  \ref{ex:pro-Omega}, we obtain a
chain map of pro-complexes
\begin{align}\label{pro-e}
\{\Omega_{S/I^m}\}_m \cong\ &
\Omega_{S}/I^\infty\Omega_S = S/I^\infty\otimes\wedge^* V = \\
S\otimes_{S^e} L/J^\infty L \ \map{\epsilon_S}\ &
S\otimes_{S^e}P_*\cong \{S/I^m\otimes_{S^e}C^{bar}_*(S/I^m)\}_m.\notag
\end{align}
where the boundary operator of the first two complexes is the zero map.  We
observe that the homology of the right side computes $HH_*(S/I^\infty)$,
and that \eqref{pro-e} induces the same map in
homology as the map \eqref{pro-hkr} (obtained from the shuffle product).

We shall prove that \eqref{pro-e} is a quasi-isomorphism.
By \ref{cor:pro-acyc}, we have a map

\[
\mu: P_*=S^e_\alpha/J^\infty\otimes_{S^e}C^{bar}_*(S/I^\infty)\to L/J^\infty L
\]
which again covers the identity of $S/I^\infty S$. Note moreover that
\[
L_0/J^\infty L_0= S^e_\alpha/J^\infty=
S^e_\alpha/J^\infty\otimes_{S^e}C_0^{bar}(S/I^\infty S)
\]
so we can take $\mu$ to be the identity in degree zero.
By the uniqueness in \ref{cor:pro-acyc}, 
we obtain an $S^e_\alpha$-linear homotopy $\mu\epsilon\to 1_{L/J^\infty L}$.
Applying $S\otimes_{S^e_\alpha}$ and taking homology in \eqref{pro-e},
we obtain that \eqref{pro-hkr} is a monomorphism in pro-homology.

To prove that \eqref{pro-hkr} is also surjective we note that, by
\eqref{equivfil}, we may augment the pro-complex
$Q_*$ 
to a pro-acyclic complex by adding $S/I^\infty S$ in degree $-1$.
Now consider  the map $\{P_{*,2m}\}_m\to \{Q_{*,m}\}$
which is zero in degrees $-1,0$ and agrees with $\epsilon\mu-1$
in higher degrees.  Applying Lemma \ref{pro-homo} to this map
(relative to the ideal $J^2$), we see that the composition with
\eqref{equivfil} is zero on homology. Since this map is $\epsilon\mu-1$,
$\epsilon\mu$ is the identity map on homology, as required.
\end{proof}
\end{subequations}

Note that the map $\mu$ constructed in the proof depends critically
upon the Artin-Rees Lemma for $\cJ$, even in degree 2.
The following example illustrates
Theorem \ref{thm:pro-hkr} in a simple case. Set
$\Lambda_m=F[x]/(x^{m+1})$ and $\Lambda_m^R=R[x]/(x^{m+1})$.
It is well known that $HH_n(\Lambda_m)\cong\Lambda_{m-1}$ for all $n>0$.
By the K\"unneth formula,
$HH_*(\Lambda^R_m)=HH_*(R)\otimes_F HH_*(\Lambda_m)$.

\begin{lem}\label{lem:TCHH}
Let $F$ be a field of characteristic zero.
For any $F$-algebra $R$ and for all $M>2m$, the image of
$\HHF_n(\Lambda^R_{M}) \to \HHF_n(\Lambda^R_m)$ is
$\bigl(\HHF_n(R)\otimes\Lambda_m\bigr) \oplus
\bigl(\HHF_{m-1}(R) \otimes \Omega^1_{\Lambda_m/F}\bigr)$.

In particular, if $R$ is 
smooth over $F$, then the image is
\[ \Omega^n_{\Lambda^R_m/F}=
\bigl(\Omega^n_{R/F}\otimes_F\Lambda_m\bigr) \oplus
\bigl(\Omega^{n-1}_R\otimes_F\Omega^1_{\Lambda_m/F}\bigr).
\]
\end{lem}

\begin{proof}
Explicit $\Lambda_m$-module generators were given in \cite[1.10]{GRW} for
$\HHF_*(\Lambda_m)$ in terms of $u=[x]$ and of the sum
$t_m$ of all terms $x^a[x^b|x]$ with $a+b=m$, considered as elements of
the bar complex: $\HHF_{2i}(\Lambda_{m})$ is generated by $xt_m^i$, and
$\HHF_{2i+1}(\Lambda_m)$ is generated by $ut_m^i$. If $M>2m$, every term in
$t_{M}$ contains an $x^{m+1}$ factor

and so vanishes in the bar complex of $\Lambda_m$. It follows that the map
$\HHF_n(\Lambda_{M})\to \HHF_n(\Lambda_m)$
is zero for all $n>1$. By the K\"unneth formula, the map from
$\HHF_n(R\otimes_F\Lambda_{M})$ to $\HHF_n(R\otimes_F\Lambda_{m})$
vanishes on all summands except for $\HHF_n(R)\otimes_F\Lambda_{M}$ and
$\HHF_{n-1}(R)\otimes_F\Omega^1_{\Lambda_M/F}$,
where it is the natural surjection.
\end{proof}

We can use the following spectral sequence to replace $F$ by any
subfield $k$ in the statement of Theorem \ref{thm:pro-hkr}.

\begin{lem}\label{lem:main11} (Kassel-Sledsj\oe, \cite[4.3a]{kasle})
Let $k\subseteq F$ be fields of characteristic zero.
For each $p\ge 1$ there is a bounded second quadrant homological
spectral sequence ($0\le i<p$, $j\ge0$):
\[ \qquad \qquad \qquad
{}_p E^1_{-i,i+j}=\Omega^{i}_{F/k}\otimes_FHH^{(p-i)}_{p-i+j}(R/F)
\Rightarrow HH_{p+j}^{(p)}(R/k)
\]
\end{lem}
\medskip

\begin{prop}\label{prop:hkrq}
Let $k\subset F$ be fields of characteristic zero, $S$ an algebra over
$F$, essentially of finite type and smooth, and $I\subset S$ an
ideal. Write $\Omega$ for $\Omega_{/k}$. For each $p\ge 0$, the
shuffle product induces an isomorphism of pro-$S$-modules
\begin{equation}\label{pro-hkr-k}
\{\Omega^p_{S/I^m}\}_m\to \{HH_p((S/I^n)/k)\}_m
\end{equation}
\end{prop}

\begin{proof}
Applying the natural, uniformly bounded spectral sequence of Lemma
\ref{lem:main11} levelwise, we obtain a bounded spectral sequence in
the category of pro-$S$-modules, which converges to the pro-Hochschild
homology over $k$ of $S/I^\infty$. It follows from Theorem
\ref{thm:pro-hkr} that the latter spectral sequence degenerates,
being zero for $j\ne0$, proving the result.
\end{proof}

\medskip
For any $\Q$-algebra $R$, let $(C(R),b,B)$ 
be the usual cyclic mixed complex of $R$; we have a canonical map
$e:(C(R),b,B)\to (\Omega_R,0,d)$
(see \cite[9.8.12]{WeibelHA94}).

\medskip
\begin{prop}\label{prop:hinfHN}
The map of mixed complexes $e:(C,b,B)\to (\Omega, 0,d)$ induces an
equivalence of fibration sequences for every $X$ in $\SchF$:

\[
\xymatrix{\bbH(X_\inf,HN)\ar[r]\ar[d]_\sim & HP(X)\ar[r]\ar[d]_\sim
    & \bbH(X_\inf,HC)[2]\ar[d]_\sim\\
\prod_{i\le 0}\bbH(X_\inf,\cO)[2i]\ar[r]
    & \prod_{i\in\Z}\bbH(X_\inf,\cO)[2i]\ar[r]
    & \prod_{i>0}\bbH(X_\inf,\cO)[2i]}
\]
\end{prop}

\begin{proof}
The top row is $\bbH(X_\inf,-)$ applied to the triangle
$HN\to HP\to HC[2]$, using Example \ref{ex:HPinf} below.
By Corollary \ref{cor:HComega}, the bottom row is $\bbH(X_\inf,-)$ applied to
the triangle $HN\Omega \to HP\Omega \to HC\Omega[2]$. The vertical
maps are the maps $e$, so the diagram commutes. It remains to show that
$e$ induces equivalences. For this, we may assume that $X$ is affine
by repeated applications of Lemma \ref{infZdescent}.

Consider the induced morphism
$e:\bbH(X_\inf,\cE)\to\bbH(X_\inf,\cE\Omega)$,
where $\cE$ is $HC$, $HP$ or $HN$.
By \ref{rem:notcrystal} and \ref{sheaf-holim}(2), 
$e$ is $\Tot_\nu$ applied to the map
$\varprojlim_m\cE(Y^{\nu+1}_m)\to \varprojlim_m(\cE\Omega)(Y^{\nu+1}_m)$,
which is a weak equivalence by \eqref{pro-hkr-k}.
\end{proof}

\begin{ex}\label{ex:hinfHN} Taking cohomology in \ref{prop:hinfHN},
we see that $H^n(X_\inf,HN)$ is:
0 for $n<0$; $H^0(X_\inf,\cO)$ for $n=0$; and
for $n>0$ it is the (finite) product
\[
H^n(X_\inf,HN) = \prod_{0\le j\le n/2} 
H^{n-2j}(X_\inf,\cO).
\]
\end{ex}

\bigskip
\section{Infinitesimal hypercohomology spectra}\label{sec:Hinf}

In this section, we rework the homological material of section
\ref{sec:infh} in the context of presheaves of spectra. We need this
generality in order to form the infinitesimal hypercohomology spectrum
for $K$-theory, introduced in section \ref{sec:Kinf}. The main result
of this section (Theorem \ref{thm:fibrant}) is that our construction
is the categorical hypercohomology spectrum, {\it i.\ e.}, global
sections of the fibrant replacement functor.

\begin{defn}\label{def:Hinf2}
Let $\cE$ be a presheaf of spectra on $X_\inf$.
We define $\bbH(X_\inf,\cE)$ to be the
homotopy limit
\[
\bbH(X_\inf,\cE)=\holim_{T\in\inf(X)}\bbHz(T,\cE).
\]
By construction \cite[XI.3.4]{BK}, there is a canonical map
$\bbH(X_\inf,\cE)\to \bbHz(X,\cE)$. Since this spectrum definition is
parallel to the homological construction in \eqref{sheafHinf},
we see that $\bbH(X_\inf,\cE)$ agrees with the homological definition
\eqref{def:Hinf1} when $\cE$ is an
Eilenberg-Mac\,Lane spectra associated to a complex of sheaves.
\end{defn}

\begin{ex}\label{ex:piHN}
For all $X$ in $\SchF$,
we see from \ref{prop:hinfHN} and \ref{ex:hinfHN} that
\[
\bbH(X_\inf,HN)\cong
\prod\nolimits_{j\ge0}\Omega^{2j}\bbH(X_\inf,\cO).
\]

In particular, if $n>0$ then $\pi_n\bbH(X_\inf,HN)=0$, and
$\pi_0\bbH(X_\inf,HN)=\bbH^0(X_\inf,\cO)$.
\end{ex}

Because $\holim$ and $\bbHz(X,-)$ preserve fibration sequences, we have:

\begin{lem}\label{lem:fibrations}
$\bbH(X_\inf,-)$ preserves fibration sequences.
\end{lem}

We say that $\cE$ is {\it nilinvariant} on $X_\zar$ if it takes
thickenings $U\harrow T$ to weak equivalences for every open $U$ in $X$.

\smallskip

\begin{ex}\label{Hinf-nilinvar}
For example, $\bbH(-_\inf,\cE)$ is nilinvariant on $X_\zar$ because for
any thickening $X\harrow X'$ the map $\inf(X')\to\inf(X)$ has a left
adjoint (the pushout of $U\hookrightarrow T$ along the unique $U\subset U'$),
so that $\holim_{\inf(X')}$ and $\holim_{\inf(X)}$ are weak equivalent.
\end{ex}

\begin{lem}\label{lem:Zdnil}
If $\cE$ is nilinvariant on $X_\zar$, and satisfies Zariski descent on
$X$, then
\[
\bbH(X_\inf,\cE)\cong \cE(X).
\]
\end{lem}

\begin{proof}
By nilinvariance, $\cE(T) \simeq \cE(U)$ for any infinitesimal
thickening $U\hookrightarrow T$ in $\inf(X)$. Therefore
$\bbH(X_\inf,\cE) \cong \holim_{U\in X_{\zar}} \cE(U)$.
By Zariski descent, this homotopy limit is
$\bbHz(X,\cE) \simeq \cE(X)$.
\end{proof}

\begin{ex}\label{ex:HPinf}
Periodic cyclic homology $HP$ is nilinvariant
(see Goodwillie \cite{Goodw1} \cite[9.9.9]{WeibelHA94}),
and satisfies Zariski descent (see \cite[2.9]{chsw}). Thus
$\bbH(X_\inf,HP)\cong HP(X).$
\end{ex}

\smallskip
Here is the spectrum analogue of Example \ref{ex:338} and
\ref{sheaf-holim}, which says that when $X=\Spec(A)$ is affine, the
groups $\bbH^*(X_\inf,\cE)$ may be computed using the thickenings
$Y^\nu_m=\Spec (S^{\otimes \nu}/I_\nu^m)$ constructed in \ref{ex:338}.

\begin{lem}\label{lem:affine_hinf}
Let $\cE$ be a presheaf of 
spectra on $\SchF$ satisfying Zariski
descent. Then for each $F$-algebra $A$ of finite type,
presented as in Example \ref{ex:338}:
\[
\bbH((\Spec A)_\inf,\cE)\simeq \Tot_\nu \cE(Y^{\nu+1}_\bu).
\]
\end{lem}

\begin{proof} Set $X=\Spec(A)$.
By definition \ref{def:Hinf2}, $\bbH(X_\inf,\cE)\simeq \holim_T \cE(T)$,
and each $\cE(T)$ is a fibrant spectrum;
our goal is to interpret the homotopy limit via $\Tot$.
Let $A=S/I_1$ be a presentation of $A$ as a quotient of a smooth
$F$-algebra. As in Example \ref{ex:338}, we form the cosimplicial tower
of algebras $S^{\otimes\nu+1}_\bu$ and the simplicial cotower
$Y^\nu_\bu$ in the category $\inf(X)$.

\goodbreak
By Lemma \ref{lem:5.1}, Bousfield-Kan cofinality \cite[XI.9.2]{BK}
and \cite[XI.4.1--4.3]{BK}, 
\[
\holim_T \cE(T)\simeq\holim_{\Delta}\holim_{\N^{op}}\cE(Y^{\nu+1}_m)
\simeq\holim_{\Delta}\cE(Y^{\nu+1}_\bu).
\]
As in \ref{sheaf-holim}, but using \cite[XI.4.4]{BK},
this is $\Tot_\nu$ of $\cE(Y^{\otimes\nu+1}_\bu)$.
\end{proof}

Now we compare the $\bbH(X_\inf,-)$ construction of Definition
\ref{def:Hinf2} with the categorical hypercohomology construction of
\cite{AKTEC} and \cite{JardineSPS}. As
one might expect, it turns out that infinitesimal hypercohomology can
be computed as global sections of a fibrant replacement in an
appropriate model structure.

Recall from \cite{JardineSHT} that there is a ``local injective''
 model structure on
presheaves of spectra on any site, and in particular on $\inf(X)$. A
map $\cE\to\cE'$ is an (infinitesimal) {\it local weak equivalence} if
it induces an isomorphism on sheaves of homotopy groups; it is a
cofibration if each $\cE(T)\to\cE'(T)$ is a cofibration of spectra
in the sense of \cite{BF}; fibrations
are defined by the right lifting property.

\begin{prop}\label{prop:lwe}
For any local weak equivalence $\cE\to\cE'$  of presheaves of spectra
on $\inf(X)$, the map $\bbH(X_\inf,\cE)\to\bbH(X_\inf,\cE')$
is a weak equivalence of spectra.
\end{prop}

\begin{proof}
By the discussion in section \ref{sec:infh}, the hypothesis means that
the stalks of $\pi_*\cE$ and $\pi_*\cE'$ are isomorphic at every
Zariski point $t\in T$, for any $U\hookrightarrow T$ in $X_\inf$.
Fixing $T$, this shows that $\cE_T\to\cE'_T$ is a Zariski local weak
equivalence in $T$, and hence that $\bbHz(T,\cE)\to\bbHz(T,\cE')$ is a
homotopy equivalence. By \cite[XI.5.6]{BK}, $\holim_T$ preserves
homotopy equivalences, so $\bbH(X_\inf,\cE)\simeq\bbH(X_\inf,\cE')$.
\end{proof}

Recall that a {\it fibrant replacement} of $\cE$ is a cofibration
$\cE \to \cE'$ with $\cE'$ fibrant which is a weak equivalence.
Because the $T$ form a covering sieve of the terminal sheaf $\ast$ on
$X_\inf$, the global sections spectrum of $\cE'$ is 
\[ \hom(\ast,\cE')=\hom({\varinjlim}_T T,\cE')=
{\varprojlim}_T \hom(T,\cE')={\varprojlim}_T\cE'(T).
\]
Here $\hom(T,\cE')$ denotes the usual spectrum Hom from
a presheaf of sets into a spectrum, and $T$ is regarded as
a (representable) presheaf of sets.

\begin{thm}\label{thm:fibrant}
Let $\cE$ be a presheaf of spectra on $X_{\inf}$. Suppose $\cE \to \cE'$ is
a fibrant replacement. Then $\bbH(X_\inf,\cE) \cong \hom(\ast,\cE').$
\end{thm}

\begin{proof}
By \ref{prop:lwe},
$\bbH(X_\inf,\cE)\simeq\bbH(X_\inf,\cE')$, so we may assume $\cE=\cE'$.

For any $U\subseteq T$ in $\inf(X)$, consider the functor
$t:\zar(T)\to \inf(X)$ sending a Zariski open $T'\subseteq T$ to the
evident thickening of $U\times_T T'$. It induces a morphism of topoi $t_*$,
from sheaves on $X_\inf$ to Zariski sheaves on $T$. 
Now $t_*$ preserves globally fibrant objects by \cite[p.119]{JardineLSS},
so we consider $\cE'_T=t_*(\cE')$.

Since globally fibrant presheaves on $T$ satisfy Zariski descent
\cite{JardineSPS},  
it follows that 
$\cE'(T)=\cE'_T(T)\simeq\bbHz(T,\cE'_T)$ and hence that
$\bbH(X_\inf,\cE')\simeq\holim_T\cE'(T)$.

Consider the covering sieve $\{T\}$ of $*$, the terminal sheaf on $X_\inf$.
By \cite[p.367 and Lemma 11]{JardineStacks}, $\hocolim T\to\ast$ is a
local weak equivalence. Combining this with the fact that
$\cE'$ is fibrant, the fact that $\cE'(T)=\hom(T,\cE')$ for $T\in\inf(X)$,
and adjointness \cite[XII.4.1]{BK}, we obtain homotopy equivalences:
\[
\hom(\ast,\cE') \simeq \hom(\hocolim{}_T T,\cE') \cong
\holim{}_T \hom(T,\cE') = \bbH(X_\inf,\cE').
\qedhere\]
\end{proof}

\begin{subremark} Fibrant presheaves are not always nilinvariant.
In particular, $\cE'(X)$ is not always equivalent to
$\bbH(X_\inf,\cE') \cong \hom(\ast,\cE')$.
\end{subremark}

\bigskip
\section{Infinitesimal $K$-theory}\label{sec:Kinf}

In this section, we apply the infinitesimal hypercohomology
construction to the algebraic $K$-theory spectrum.

\begin{defn}\label{def:KIT}
For $U\subset T$ in $\inf(X)$, let $KI(T)$ and $HNI(T)$ denote the
respective homotopy fibers of $K(T)\to K(U)$ and $HN(T)\to HN(U)$.
The relative Chern character induces a natural map $KI(T)\to HNI(T)$.
Goodwillie's theorem \cite{G86} states that
$KI(T)\to HNI(T)$ is a homotopy equivalence. (It is well known
that the $\pi_nKI(T)$ are uniquely divisible; see
\cite{WeibelNil}.)
\end{defn}

\begin{subremark}\label{Hinf-KOI}
Let $K(\cO,I)$ and $HN(\cO,I)$ denote the presheaves sending $U\subset T$
to the fibers of $K(\cO(T))\to K(\cO(U))$ and $HN(\cO(T))\to HN(\cO(U))$.
Then $K(\cO,I)\to KI$ and $HN(\cO,I)\to HNI$ are local weak equivalences
on $\inf(X)$, because every thickening is locally affine.
By \ref{prop:lwe}, $\bbH(X_\inf,K(\cO,I))\simeq \bbH(X_\inf,KI)$ and
$\bbH(X_\inf,HN(\cO,I)) \simeq \bbH(X_\inf,HNI)$.
\end{subremark}

\begin{thm}\label{chsquare}
The map $\bbH(X_\inf,KI)\to\bbH(X_\inf,HNI)$ is a homotopy equivalence
for all $X$ in $\SchF$,
and there is a commutative diagram whose rows are fibrations:
\begin{equation*}
\xymatrix{\bbH(X_\inf,KI)\ar[d]_\simeq\ar[r]&
\bbH(X_\inf,K)\ar[d]_{\bbH(ch)}\ar[r]&K(X)\ar[d]_{ch} \\
\bbH(X_\inf,HNI) \ar[r]& \bbH(X_\inf,HN)\ar[r] & HN(X).}
\end{equation*}
\end{thm}

\begin{proof}
Applying $\bbH(X_\inf,-)$ to $KI(T)\map{\simeq} HNI(T)$
yields the first assertion.

By Lemma \ref{lem:Zdnil}, $\bbH(X_\inf,T\!\mapsto{\!}K(U))\cong K(X)$,
and similarly for $T\mapsto HN(U)$.
By Lemma \ref{lem:fibrations},
applying $\bbH(X_\inf,\!-)$ to $KI(T)\to K(T)\to K(U)$ and its $HN$ analogue
yields the desired equivalence of fibration sequences.
\end{proof}

\begin{thm}\label{thm:unfishy}
If $n\ge1$ then $HN_{n}(X)\cong \pi_{n-1}\bbH(X_\inf,KI)$
for all $X$ in $\SchF$.
\end{thm}

\begin{proof}
The result is immediate from \ref{chsquare} when $n>1$, since
$\pi_n\bbH(X_\inf,HN)=0$ for $n\ge1$ by \ref{ex:piHN}.
For $n=1$ it suffices to show that $\pi_0\bbH(X_\inf,HN)\to HN_0(X)$
is an injection. In fact,
\[ \pi_n\bbH(X_\inf,HN)\to HN_n(X)\to HP_n(X) \]
is a split injection for all $n$, by Proposition \ref{prop:hinfHN}.
\end{proof}

\begin{subremark}\label{rem:negK}
When $n\le0$, 
$\pi_{n-1}\bbH(X_\inf,KI) \cong \prod_{i>0} \bbH^{2i-n}(X_\inf,\cO)$
for affine $X$. This follows from \ref{chsquare}, \ref{prop:hinfHN} and
$HN_n(X)=HP_n(X)\cong \prod_{i} \bbH^{2i-n}(X_\inf,\cO)$.
\end{subremark}

\smallskip
\begin{subremark}
Although the setup is analogous to that of \cite{CortInf} and
\cite{CortKABI}, those articles allow noncommutative thickenings,
and \cite{CortInf} uses connective $K$-theory.
These differences account for the restriction $n\geq 1$ in
Theorem \ref{thm:unfishy}.
Nevertheless, up to these slight differences in definitions,
Theorem \ref{thm:unfishy} recovers the second main result,
Theorem 6.2(i), of \cite{CortInf}. That is, for $n\ge1$ the Chern character
$K_n(X)\to HN_n(X)$ may be indentified with the map
$K_n(X)\to \pi_{n-1}\bbH(X_\inf,KI)$ in \ref{chsquare}.
This fails for $n=0$, since second map omits the ``rank'' component
\[ch_0:K_0(X)\to H^0(X,\Z)\to H^0(X_\inf,\cO)\]
of the Chern character $K_0(X)\to HN_0(X)$, as Remark \ref{rem:negK} shows.
\end{subremark}

\medskip
We conclude this section with a comparison (Theorem \ref{thm:hinfK})
between infinitesimal hypercohomology $\bbH(X_\inf,K)$ and the
``infinitesimal $K$-theory'' used in our earlier papers \cite{chsw}
and \cite{chw} and based upon the construction in the eponymous paper
\cite{CortInf}.
\goodbreak

\begin{defn}\label{def:Kinf}
Let $K^\inf(X)$ denote the homotopy fiber of the Chern character
$ch:K(X)\to HN(X)$. Here $K(X)$ is non-connective $K$-theory, not
the connected and rational version used in \cite{CortInf}.
\end{defn}

\begin{lem}\label{hkinf=kinf}
$\bbH(X_\inf,K^\inf)\simeq K^\inf(X).$
\end{lem}

This lemma is immediate from \ref{chsquare} and \ref{lem:fibrations}.
Alternatively, note that $K^\inf$
satisfies Zariski descent, because both $K$ and
$HN$ do.  Moreover, $K^\inf$ is nilinvariant, by \ref{def:KIT},
so \ref{hkinf=kinf} also follows from Lemma \ref{lem:Zdnil}.

\begin{thm}\label{thm:hinfK}
For any $X\in\SchF$, there is a homotopy fibration sequence:
\begin{equation*}\label{hik-kinf}
K^\inf(X)\to \bbH(X_\inf,K)\to
\prod\nolimits_{i\ge0}\Omega^{2i}\bbH(X_\inf,\cO)
\end{equation*}
In particular, for $n\ge 1$, $K_n^\inf(X)\cong\pi_n\bbH(X_\inf,K)$.
\end{thm}

\begin{proof}
Applying $\bbH(X_\inf,-)$ to the fibration of \ref{def:Kinf}, and using
\ref{hkinf=kinf} and \ref{lem:fibrations}, we get a fibration sequence
$K^\inf(X)\to \bbH(X_\inf,K)\to\bbH(X_\inf,HN)$. The sequence in
\ref{hik-kinf} and the homotopy group calculation follow from Example
\ref{ex:piHN}.
\end{proof}

\begin{subremark}
One should compare Theorem \ref{thm:hinfK} to
\cite[Theorem 6.2(ii)]{CortInf}. This is not straightforward, because
our notation is not compatible with that in \cite{CortInf}.
\end{subremark}

\section{Cathelineau's Theorem for Schemes}\label{sec:Cath-schemes}

In this section, we develop both a scheme-theoretic version
(\ref{KUIcompatible}) and an infinitesimal version (\ref{Hinf-KI}) of
Cathelineau's Theorem (see Theorems \ref{thm:cath} and \ref{thm:KOIcompatible}
in the Appendices). Although we do not need the scheme-theoretic version
for our Main Theorem \ref{main-intro}, it is of independent interest.
It also sets the stage for the infinitesimal version,
which we will need for Theorem \ref{thm:mth}, and hence for our Main Theorem.
Throughout this section, we work over a field of characteristic~0.

Let $I$ be a nilpotent sheaf of ideals on a scheme, and consider
the Zariski presheaf $K(-,I)$, sending $U$ to the homotopy fiber of
$K(U)\to K(U/I)$, where $U/I$ is the scheme over $U$ with structure sheaf
$\cO_U/I$.  By \cite{TT}, $K$ and $K(-,I)$ satisfy Zariski descent.
Since $K(\cO,I)\to K(-,I)$ is a local weak equivalence for the Zariski
topology, it follows that we have $K(U,I)\simeq\bbH_\zar(U,K(\cO,I))$.

By \ref{deloop-lambda}, 
both $\lambda^k$ and $\psi^k$ are spectrum maps from $K(\cO,I)$ to
itself. Applying Zariski descent, they induce spectrum maps from each
$K(U,I)$ to itself.

We define $HC(-,I)$ and $HN(-,I)$ similarly, or as the Eilenberg-Mac\,Lane
spectra associated to the appropriate chain complexes of sheaves.
Since $\lambda^k$ and $\psi^k$ are associated to chain maps,
compatible with the homotopy equivalences $B:HC(U,I)[1]\to HN(U,I)$,
they are spectrum maps compatible with $B$.

\begin{thm}\label{KUIcompatible}
If $I$ is a nilpotent sheaf of ideals on a scheme,
then the relative Chern character $K(-,I) \to HC(-,I)[1]\simeq HN(-,I)$
is compatible with the operations $\lambda^k$ in the sense that
$\psi^k(ch\, x) = k\cdot ch(\psi^k x)$ for each $x\in K_m(U,I)$.

In addition, we have a homotopy commutative diagram of presheaves of spectra:
\begin{equation*}
\xymatrix{ K(-,I) \ar[r]^{\simeq} \ar[d]^{\lambda^k} &
	 HC(-,I)[1] \ar[d]^{\lambda^k} \ar[r]^{\simeq}_B  &
	HN(-,I) \ar[d]^{\lambda^k} \\
K(-,I) \ar[r]^{\simeq} & HC(-,I)[1]\ar[r]^{\simeq}_B & HN(-,I).
}\end{equation*}
\end{thm}

\begin{proof} By the spectrum version \ref{deloop-lambda} of
Theorem \ref{thm:KOIcompatible}, we have a homotopy commutative diagram
of presheaves of spectra:
\begin{equation*}
\xymatrix{ K(\cO,I) \ar[r]^{\simeq} \ar[d]^{\lambda^k}&
	HC(\cO,I)[1] \ar[r]^{\simeq}_B  \ar[d]^{\lambda^k} &
	HN(\cO,I) \ar[d]^{\lambda^k} \\
K(\cO,I) \ar[r]^{\simeq} & HC(\cO,I)[1] \ar[r]^{\simeq}_B & HN(\cO,I).
}\end{equation*}
Applying Zariski hypercohomology (and homotopy groups), we get the result.
\end{proof}

\smallskip\noindent
Let $K^{(i)}(-,I)$, $HC^{(i-1)}(-,I)$ and $HN^{(i-1)}(-,I)$ denote the
respective homotopy fibers of $\psi^k-k^i$ on $K(-,I)$, $HC(-,I)$ and
$HN(-,I)$. 
Copying the proof of Corollary \ref{prodKiOI} proves the following.

\begin{cor}\label{KiUIcompatible}
There is a homotopy commutative diagram for each $U$:
\begin{equation*}
\minCDarrowwidth10pt
\begin{CD}
K(U,I) @>{\simeq}>> HC(U,I)[1] @>{\simeq}>{B}> HN(U,I) \\ \vspace{-2pt}
@VV{\simeq}V @VV{\simeq}V  @VV{\simeq}V  \\
\prod_{i=1}^\infty K^{(i)}(U,I)  @>{\simeq}>>
 \prod_{i=1}^\infty HC^{(i-1)}(U,I)[1]  @>{\simeq}>{B}>
 \prod_{i=1}^\infty HN^{(i)}(U,I).
\end{CD}\end{equation*}
\end{cor}

\medskip
We now turn to the analogues of Theorem \ref{KUIcompatible} and
Corollary \ref{KiUIcompatible} for the infinitesimal topology on
a fixed $X$ in $\SchF$. In order to distinguish the tautological ideal
of a thickening $U\subset T$ from the Zariski case, we adopt the
notation $KI$ and $HNI$ from \ref{def:KIT}
rather than the notation $K(-,I)$ and $HN(-,I)$ used above.

The relative Chern character $ch:KI\to HNI$ induces an infinitesimal
Chern character $\bbH(X_\inf,KI) \to \bbH(X_\inf,HNI)$, and we saw in
Theorem \ref{chsquare} that it is a homotopy equivalence.
Let $KI^{(i)}(T)$ and $HNI^{(i)}(T)$ denote the respective
homotopy fibers of $\psi^k-k^i$ on $KI$ and $HNI$;
as in \ref{KiOI} these are defined using a fixed $k\ge2$, but are independent
of this choice up to homotopy equivalence.

\begin{thm}\label{Hinf-KI}
For any $X$ in $\SchF$, the infinitesimal Chern character
is compatible with the operations $\lambda^k$. That is,
we have a commutative diagram of spectra:
\begin{equation*}
\xymatrix{ \bbH(X_\inf,KI) \ar[r]_{ch}^{\simeq} \ar[d]^{\lambda^k} &
	\bbH(X_\inf,HNI) \ar[d]^{\lambda^k} \\
\bbH(X_\inf,KI) \ar[r]^{\simeq}_{ch} &
 \bbH(X_\inf,HNI).
}\end{equation*}
\noindent
There is also a homotopy commutative diagram:
\begin{equation*}
\minCDarrowwidth10pt
\begin{CD}
\bbH(X_\inf,KI) @>{\simeq}>{ch}>  \bbH(X_\inf, HNI) \\ \vspace{-2pt}
 @VV{\simeq}V @VV{\simeq}V \\
\prod_{i=1}^\infty \bbH(X_\inf,KI^{(i)})  @>{\simeq}>{ch}>
\prod_{i=1}^\infty \bbH(X_\inf,HNI^{(i)}).
\end{CD}\end{equation*}
\end{thm}

\begin{proof}
The maps $K(\cO,I)\to KI$ and $HN(\cO,I)\to HNI$ are local weak equivalences
on $X_\inf$ by \ref{Hinf-KOI}, and are compatible with the operations
$\psi^k$.  Hence they induce local weak equivalences
$K^{(i)}(\cO,I)\to KI^{(i)}$ and $HN^{(i)}(\cO,I)\to HNI^{(i)}$.
By Proposition \ref{prop:lwe}, 
the diagrams in question are weak equivalent to the
functorial infinitesimal hypercohomology construction \ref{def:Hinf2}
applied to the commutative diagrams in Theorem \ref{thm:KOIcompatible},
Corollary \ref{deloop-lambda} and Corollary \ref{prodKiOI}.
\end{proof}


\section{Hodge decomposition}\label{sec:main}

We now turn to the proof of the main theorem \ref{main-intro}
stated in the introduction. We first prove the result for $n\ge1$,
and for all $X$ in $\SchF$, where $F$ is a field of characteristic~0.

\begin{thm}\label{thm:mth}
For all $X$ in $\SchF$ and all $n\ge1$,
the map $ch:K_n(X)\to HN_n(X)$ satisfies $\psi^k(ch) = k\cdot ch(\psi^k)$,
i.e., $ch$ sends $K_n^{(i)}(X)$ to $HN_n^{(i)}(X)$.
\end{thm}

\begin{proof}
The left half of the following diagram commutes by Theorem \ref{chsquare},
and the right half of this diagram commutes by Theorem \ref{Hinf-KI}.
The lower left map is injective when $n\ge1$ because in that case
$\pi_n\bbH(X_\inf,HN)=0$ by \ref{ex:piHN}.
\begin{equation*}
\minCDarrowwidth10pt
\begin{CD}
K_n(X) @>{\partial}>> \pi_{n-1}\bbH(X_\inf,KI) @>{\cong}>>
	\prod_{i=1}^\infty \pi_{n-1}\bbH(X_\inf,KI^{(i)})\\
@VV{ch}V @VV{\cong}V @VV{\cong}V \\
HN_n(X) @>{\partial}>{\text{into}}> \pi_{n-1}\bbH(X_\inf, HNI)
	@>{\cong}>> \prod_{i=1}^\infty \pi_{n-1}\bbH(X_\inf,HNI^{(i)}).
\end{CD}\end{equation*}
The two horizontal maps $\partial$ commute with $\psi^k$ because
the operations $\psi^k$ are compatible with the fibration
$KI(T)\to K(T)\to K(U)$ (and its $HN$ analogue) used in the proof of
Theorem \ref{chsquare}. Thus if $x\in K_n^{(i)}(X)$ then
$\partial(x)\in \pi_{n-1}\bbH(X_\inf,KI^{(i)})$, and a diagram chase
shows that $ch(x)$ is in $HN_n^{(i)}(X)$.
\end{proof}

\begin{cor}\label{cor:ch and psik}
For all $X$ in $\SchF$ and all $n$,
the map $ch:K_n(X)\to HN_n(X)$ satisfies $\psi^k(ch) = k\cdot ch(\psi^k)$,
i.e., $ch$ sends $K_n^{(i)}(X)$ to $HN_n^{(i)}(X)$.
\end{cor}

\begin{proof} We proceed by downward induction on $n$,
Theorem \ref{thm:mth} being the base case $n\ge1$.
Given $x\in K_n(X)$, the element $\{x,t\}$ of $K_{n+1}(X\times\Gm)$
satisfies $\psi^kch(\{x,t\}) = k\cdot ch(\psi^k\{x,t\})$ by induction.
But $\psi^k$ is multiplicative on $K$-theory,
and the Chern character is multiplicative by \cite[\S5]{HJ} \cite{Ginot04},
so by \ref{ex:ch(t)} we have:
\[
ch(\psi^k\{x,t\}) = ch(\psi^kx)\cdot ch(\psi^kt)
 = k\, ch(\psi^kx)\cdot ch(t).
\]
Although the Adams operations $\psi^k$ are not multiplicative on $HN$,
the operations $\bar{\lambda}^k=(-1)^{k-1}\lambda^k$ are
multiplicative because they commute with the shuffle product (see
\cite[4.5.14, 5.1.14]{LodayHC92}).  Since the Adams operations satisfy
$\psi^k=k\bar\lambda^k$, we have $k\psi^k(a\cdot
b)=\psi^ka\cdot\psi^kb$ in $HN$. Hence
\[
k\psi^k(ch\{x,t\}) = k\psi^k\bigl(ch(x)\cdot ch(t)\bigr)) =
\psi^kch(x)\cdot\psi^k ch(t)  = k^2 \psi^k ch(x)\cdot ch(t).
\]
The result follows, because we can divide by $k$ in $HN_*$,
and multiplication by
$ch(t)$ is an injection from $HN_n(X)$ into $HN_{n+1}(X\times\Gm)$.
\end{proof}

\begin{proof}[Proof of Theorem \ref{main-intro}]
Every element of $K_n(A)$ comes from $K_n(A_0)$ for some subalgebra $A_0$
of finite type, so we may assume that $A$ is of finite type over a field $F$.
If $n\le 0$, the theorem follows from \cite{Weibel93}, and the case $n\ge1$
is handled by Theorem \ref{thm:mth} above (with $X=\Spec A$).
\end{proof}

\newpage
\addtocounter{section}{-7}
\renewcommand{\thesection}{\Alph{section}} 
\section{Appendix: Cathelineau's Theorem}\label{app:A}

In this appendix, we correct the proof of Cathelineau's theorem \cite{C91}.
Recall that Goodwillie's theorem \cite{G86} identifies the
relative $K$-theory  $K_n(A,I)$ and cyclic homology $HC_{n-1}(A,I)$
of a nilpotent ideal $I$ in a $\Q$-algebra $A$. Also recall that these
groups are the direct sum of their $k^i$-eigenspaces
$K^{(i)}_n(A,I)$ and $HC^{(i)}_n(A,I)$ for the Adams operation $\psi^k$.

Actually, Goodwillie gave two such isomorphisms, the relative Chern
character $ch_*$ and the rational homotopy character $\rho_*$. These group
isomorphisms ere shown to be identical in \cite{cw-agree}, since the
maps which induce them are naturally homotopic.

\begin{Cthm}\label{thm:cath}
Let $I$ be a nilpotent ideal in a commutative $\Q$-algebra $A$.
Then Goodwillie's isomorphism $K_n(A,I) \cong HC_{n-1}(A,I)$
is an isomorphism of trivial $\gamma$-rings. That is, it
identifies $K^{(i)}_n(A,I)$ and $HC^{(i-1)}_{n-1}(A,I)$.
\end{Cthm}

Using the Science Citation Index to follow Ariadne's Thread, we see that
Cathelineau's Theorem \ref{thm:cath} has been used in
\cite{GW94}, \cite{VP95}, \cite{EV96, EV99}, \cite{Ginot04}, \cite{Kwz99}
and the present paper. In addition, \cite{LodayHC92} and
several other papers have cited it
without using it: \cite{Weibel93}, \cite{C96, C98}, \cite{Gonch1, Gonch2},
\cite{L03} and \cite{Krishna2, Krishna3}.  

The problem with the proof in \cite{C91} is that it references the
unpublished Ogle-Weibel preprint \cite{OW}, which mistakenly asserts
in \cite[1.6]{OW} that the $GL(\Q)$-invariant Chevalley-Eilenberg
chain complex $x^{GL(\Q)}(A,I)$ is quasi-isomorphic to a certain
homotopy colimit. To fix it, we shall systematically use the
combinatorial version $X(A,I)$ corresponding to the
symmetric group $\Sigma_\infty$.

\subsection*{Step 1}
Most of this step is presented in \cite{LodayHC92}.  If $A$ is a ring
with unit, and $T_n(A)$ is the subgroup of $GL_n(A)$ consisting of
upper triangular matrices, then the {\it Volodin space} $X(A)\subset
BGL(A)$ is the union of the spaces $X_n(A)=\bigcup_{\sigma\in\Sigma_n}
BT_n^\sigma(A)$; Suslin proved in \cite{Su81} that $X(A)^+$ is
contractible and that there is a homotopy fibration
\begin{equation}\label{eq:Vol}
X(A) \to BGL(A) \to BGL(A)^+.
\end{equation}

If $I$ is an ideal in $A$, we may consider the subgroup $T_n(A,I)$ of
$GL_n(A)$ consisting of matrices which are upper-triangular
modulo $I$, i.e., the preimage of the upper triangular matrices $T_n(A/I)$.
By definition, $X_n(A,I)$ is the union
$\bigcup_{\sigma\in\Sigma_n} BT_n^\sigma(A,I)$.

Let $\overline{GL}(A/I)$ denote the image of $GL(A)$ in $GL(A/I)$,
and define $K(A,I)$ to be the homotopy fiber of
$BGL(A)^+\to B\overline{GL}(A/I)^+$;
by construction, $K(A,I)$ is a connected space whose homotopy groups
$\pi_nK(A,I)$ are the relative $K$-groups $K_n(A,I)$ for all $n\ge1$.
The following theorem is Theorem 6.1 of \cite{OW}; it is cited in \cite{C91}
as well as in \cite{AO}. As usual,
$X(A,I)$ denotes the union of the $X_n(A,I)$.

\begin{thm}\label{OW6.1}
If $I$ is an ideal in $A$, there are homotopy fibrations
\begin{gather*} X(A,I) \to BGL(A) \to B\overline{GL}(A/I)^+, \\
X(A) \to X(A,I) \to K(A,I).
\end{gather*}
Moreover, there is a homology isomorphism $X(A,I)\to K(A,I)$ and a
homotopy equivalence $X(A,I)^+ \simeq K(A,I)$.
\end{thm}

\begin{proof} The original proof of Theorem 6.1 in \cite{OW} works, and
is reproduced on page 361 of \cite{LodayHC92}.
\end{proof}

\begin{lem}\label{lem:embed}
If $\sigma\in\Sigma_{n+1}$, the canonical embedding $i_n:X_n\to X_{n+1}$
and the embedding $\sigma i_n\sigma^{-1}$ are naturally homotopic.
Hence $\Sigma_\infty$ acts trivially on $H_*X(A,I)$.
\end{lem}

\begin{proof} As noted in \cite[3.1]{OW} \cite[I.5]{G85},
Suslin's argument in \cite[1.5]{Su83} goes through.
\end{proof}

\begin{cor}\label{OW6.1.1}
The direct sum of matrices makes $X(A,I)^+$ into an $H$-space.
In particular, $H_*(X(A,I);\Q)$ is a Hopf algebra.
\end{cor}

\begin{proof} Lemma \ref{lem:embed} is the ingredient needed for the
classical proof to go through. (Cf.\ \cite[Remark 6.1.1]{OW}
and the proof in \cite[11.3.4]{LodayHC92}.)
\end{proof}

In order to view $X(A,I)$ as a homotopy colimit, we need to expand
the indexing set from $\Sigma_n$ to the collection $\cA_n$
of all partial orders $\sigma$ of $\{1,...,n\}$, so that family of
subgroups $T^{\sigma}$ is closed under intersection. This expands
$\Sigma_\infty=\cup\Sigma_n$ to $\cA=\cup\cA_n$, and is the
standard adjustment described in \cite[p.~386]{G85}, \cite[1.5]{OW}
and \cite[1.2]{AO}. The point is that the natural map of simplicial
chain complexes
\begin{equation}\label{BThocolim}
\hocolim_{\cA_n} C_*(BT_n^\sigma(A,I)) \to C_*(X_n(A,I))
\end{equation}
is a quasi-isomorphism for all $n$, including $n=\infty$.
The proof of Claim III.9 in \cite{G85}
goes through in this setting, as does the proof of \cite[1.4]{OW}.

\subsection*{Step 2}
On pp.~391--393 of \cite{G85}, Goodwillie considered the Lie algebra
$\frakt_n(A)$ of strictly upper triangular matrices and introduced the sum
$x_n(A)=\sum_{\sigma\in\Sigma_n} C_*\frakt_n^\sigma(A)$
of the Chevalley-Eilenberg chain complexes as a useful homological tool.
He also introduced the subgebra $\frakt_n^\sigma(A)$ for each partial
order $\sigma$ of $\{1,...,n\}$ and proved (Claim III.9 in \cite{G85})
that the canonical map
$\hocolim_{\cA_n}C_*\frakt_n^\sigma(A) \to x_n(A)$ is a quasi-isomorphism.

There are quasi-isomorphisms
$C_*(T_n^{\sigma}(A),\Q)\map{\simeq} C_*\frakt_n^\sigma(A)$,
arising from rational homotopy theory,
which are natural in $n$ and $\sigma$. They were used by Goodwillie
in \cite[p.~392]{G85} and shown to be natural in \cite[5.11]{SW}.
As observed on p.~85 of \cite{SW}, naturality
implies that if $\Q\subseteq A$ there is a natural quasi-isomorphism
\begin{equation}\label{eq:hocolim}
C_*(X_n(A),\Q) 
\map{\simeq} x_n(A).
\end{equation}

If $I$ is a nilpotent ideal in $A$, and $\sigma$ is a partial ordering
of $\{1,...,n\}$, we may consider the nilpotent Lie subalgebra
$\frakt_n^\sigma(A,I)$ of $\mathfrak{gl}_n(A)$ consisting of matrices
whose reductions modulo $I$ are in $\frakt_n^\sigma(A/I)$. The
rational homotopy theoretic argument given by Goodwillie on pp.~392--3
of \cite{G85} (and made natural in \cite[5.11]{SW})
works for any nilpotent Lie algebra over $\Q$, and so works in
this context to yield canonical quasi-isomorphisms
\begin{equation}\label{eq:RHT}
 C_*(T_n^\sigma(A,I),\Q) \map{\simeq}
C_*\frakt_n^\sigma(A,I). 
\end{equation}

Now consider the chain complex
$x_n(A,I)=\sum_{\sigma} C_*(\frakt_n^\sigma(A,I))$.
The following lemma repairs the mistake in \cite[1.6]{OW}, which is
also a gap in the sketch of \cite[11.3.15]{LodayHC92}.

\begin{lem}\label{lem:OW1.6}
The map $\hocolim_{\cA_n} C_*\frakt_n^\sigma(A,I) \to x_n(A,I)$
is a quasi-isomorphism.
\end{lem}

\begin{proof}
Choosing a basis for $I$ and completing it to a basis of $A$,
we get a basis of $C_*\mathfrak{gl}_n(A)$ which restricts to a basis
$B_*(\sigma)$ of each $C_*\frakt_n^\sigma(A,I)$, and this family of
chain complexes 
is closed under intersection. Thus for each $p$ the functor
$\sigma\mapsto C_p\frakt_n^\sigma(A,I)$ is the free vector space
on the underlying basis functor $B_p$ from $\cA_n$ to sets. For each basis
element $b\in\cup B_p(\sigma)$ there is a unique minimal partial order
$\alpha$ so that $b\in B_p(\alpha)$.  This is enough
for the proof of Claim III.9 in \cite{G85} to go through.
\end{proof}
\goodbreak

\begin{rem} The construction of the basis $B_p(\sigma)$ is also enough
for the following proof to work; it is the original proof of \cite[1.4]{OW}.
For each $q$, the simplicial homotopy
\[
h_j(b,\sigma_0\to\cdots)=(b,\alpha=\cdots\alpha\to\sigma_j\to\cdots)
\]
from the identity of $\hocolim C_q$ to the retraction onto the subgroup
$x_q(A,I)$ shows that $H_p(\hocolim C_q)$ is: $x_q(A,I)$ if $p=0$ and zero
otherwise. Thus the spectral sequence
$E_{pq}^1=H_p(\hocolim C_q)\Rightarrow H_{p+q}\hocolim C_*$
degenerates at $E^2$ to yield the conclusion
$H_px_n(A,I)\cong H_p\hocolim C_*$ of Lemma \ref{lem:OW1.6}.
\end{rem}

\begin{thm}\label{OW2.4} There is a natural isomorphism
$H_*(X_n(A,I),\Q) \cong H_*x_n(A,I)$ for each $n$,
induced by $\Sigma_n$-equivariant quasi-isomorphisms:
\[
C_*(X_n(A,I),\Q) \simeq \hocolim_{\cA_n} C_*(T_n^\sigma(A,I),\Q) \simeq
\hocolim_{\cA_n} C_*\frakt_n^\sigma(A,I) \simeq x_n(A,I).
\]
\end{thm}

\begin{proof} Combine \eqref{BThocolim} and the homotopy colimit of
\eqref{eq:RHT} with Lemma \ref{lem:OW1.6}.
\end{proof}

\begin{rem}\label{OW2.4bis}
The natural maps in \eqref{eq:RHT} assemble directly to yield the
natural quasi-isomorphism $\phi:C_*(X_n(A,I),\Q) \map{\simeq} x_n(A,I)$.
This follows from repeated use of Mayer-Vietoris sequences,
as observed in \cite[p.~85]{SW}. A second $\hocolim$-free proof of
Theorem \ref{OW2.4} is sketched in \cite[11.3.15]{LodayHC92}.
\end{rem}

The composition $\rho$ of the Hurewicz map $K(A,I)\map{h} C_*K(A,I)$,
the homotopy equivalence of \ref{OW6.1}, the map $\phi$ of \ref{OW2.4bis}
and the Loday-Quillen map $\theta:x(A,I)\!\to\!HC(A,I)[1]$ (defined in
\cite[10.2.3/11.3.12]{LodayHC92}) and $B:HC(A,I)[1]\simeq HN(A,I)$ is:
\begin{equation}\label{eq:KAI-HC}
K(A,I) \smap{h}
C_*K(A,I) \, {\buildrel\simeq\over\leftarrow}\, C_*X(A,I)
 \smap{\phi} x(A,I) \smap{\theta} HC(A,I)[1] \smap{B} HN(A,I).
\end{equation}
\noindent Here the cyclic homology complex is taken over $k=\Q$, and
the non-negative chain complexes are regarded as simplicial sets using
the Dold-Kan correspondence.  
We write $ch^-_\rht$ for the composition
$B\theta\circ\phi: C_*X(A,I)\to HN(A,I)$
in \eqref{eq:KAI-HC}; it is natural in $A$ and $I$.
\smallskip

\begin{defn}\label{def:relChern}
(\cite{G86}; cf.\ \cite[11.3.1]{LodayHC92})
Let $I$ be a  nilpotent ideal in a $\Q$-algebra $A$.
The map $\rho:K(A,I)\to HC(A,I)[1]$ of \eqref{eq:KAI-HC} is called the
{\it rational homotopy theory character}. 
Goodwillie proved in \cite{G86} that $ch^-_\rht$ and hence $\rho$ are
homotopy equivalences for all $(A,I)$.
It is the map invoked by Cathelineau in \cite[p.600]{C91}.
\end{defn}

Here is the correction to Proposition~4.3 of \cite{OW}, which is also
the conclusion (1.2.4) of Aboughazi-Ogle \cite{AO}. Since their proof
formally relies on some assertions in \cite{OW}, this theorem
also fills in the details of the presentation in \cite{AO}.

\begin{thm}\label{OW4.3} (Cf.\ \cite[11.3]{LodayHC92})
If $I$ is a nilpotent ideal in a $\Q$-algebra $A$ then the map $\rho$
induces isomorphisms for all $m\ge1$:
\[ K_m(A,I) \cong \Prim H_m (C_*X(A,I),\Q)
\cong \Prim H_m x(A,I) \cong HC_{m-1}(A,I). \]
\end{thm}

\begin{proof}
We need only verify that the proof of (1.2.4) in \cite{AO} works {\it
mutatis mutandis}.  Equation (1.2.1) in \cite{AO} is the combination
of \eqref{BThocolim} and Theorem \ref{OW6.1}.  Equation (1.2.2) in
\cite{AO} is just Theorem \ref{OW2.4}.  Given these substitions and
Lemma \ref{lem:embed}, the proof of Claim~1.2.3 in \cite{AO} is valid.
Therefore Theorem~1.1.12 of \cite{AO} applies to the complex $x(A,I)$,
as asserted in \cite{AO}, to yield the desired isomorphisms.
\end{proof}

We will need the following result, proven in \cite[6.5]{cw-agree}.
The {\it relative Chern character} $ch: K(A,I)\to HN(A,I)$ was defined
by Goodwillie in \cite[II.3.3]{G86} as the composition of the map
$K(A,I) \map{h} C_*K(A,I) \ {\buildrel\sim\over\leftarrow} \ C_*X(A,I)$
of \eqref{eq:KAI-HC} and a natural map
$ch^- : C_*X(A,I) \to HN(A,I)$
(see \cite[11.4.6]{LodayHC92} for the definition of $ch^-$).
\goodbreak

\begin{prop}\label{ch=rht} (\cite[6.5.1]{cw-agree})
Let $I$ be a nilpotent ideal in a $\Q$-algebra $A$.
Then the map $ch^-$ 
is naturally chain homotopy equivalent to $ch^-_\rht$.
Hence the relative Chern character $ch$ and the rational homotopy
character $\rho$ (composed with $HC(A,I)[1] \simeq HN(A,I)$) are
homotopic for each $(A,I)$.
\end{prop}

\smallskip
\subsection*{Step 3}
Cathelineau's paper \cite{C91} refers to the preprint \cite{OW}
in the following places.

\begin{enumerate} 
\item On p.~597, Cathelineau cites \cite[2.3]{OW} for the
isomorphism $H_*\frakt_n(A,I)\cong H_*(T_n(A,I),\Q)$.
This isomorphism is rederived in \eqref{eq:RHT} above, and in
\cite[11.3.14]{LodayHC92}.

\item Set $G=GL(\Q)$. On p.~599, the space
$X^{G}(A,I)=\bigcup_{g\in G}BT(A,I)^g$ and the chain complex
$x^G(A,I)=\sum_{g\in G} C_*(\mathfrak{t}^g)$ of
\cite{OW} are introduced, and then \cite[2.4]{OW}
is cited for the existence of a $G$-equivariant isomorphism
$\Phi^G:H_*(x^G(A,I))\to H_*(X^G(A,I))$ in the proof of Lemma~2.

The existence of $\Phi^G$ is unclear.
As explained above, we need to replace $X^G(A,I)$ by the space
$X(A,I)$, and $x^G(A,I)$ by the chain complex $x(A,I)$. By appealing
to Theorem \ref{OW2.4} above, we get the required isomorphism
$H_*(x(A,I)\cong H_*(X(A,I))$. This corrects the proof of
Lemma 2 of \cite{C91}.

\item In the definition of $\lambda^k$ (\cite{C91}, bottom of p.~599),
Cathelineau cites \cite[3.2]{OW} for the fact that $G$ acts trivially on
the homology of $X^G(A,I)$ and $x^G(A,I)$, so that the maps $\lambda^k_n$
are compatible with inductive limits. This is not a problem any more,
by Lemma \ref{lem:embed} and \ref{OW2.4}.

\item He cites \cite[6.1]{OW} for the fact that $X^G(A,I)\to K(A,I)$
is a homology isomorphism. This is replaced by Theorem \ref{OW6.1} above.

\item In the middle of p.~600, he cites \cite[4.3]{OW} that the primitives of
$H_*(X)$ agree with $HC_{*-1}(A,I)$. This is replaced by Theorem \ref{OW4.3}.

\item The $H$-space structure on $X(A,I)^+$, used to define the topological
maps $\lambda^k_n$ from $X_n(A,I)$ to $X(A,I)^+$ on p.~601, is deduced from
Remark 6.1.1 of \cite{OW}. This is addressed in \ref{OW6.1.1} above.

\item In order to see that $\lambda^k_{n+1}$ is freely homotopic to
$i_n\circ\lambda^k_n$, Cathelineau invokes the argument of \cite[3.1]{OW}.
In fact, these maps differ by conjugation by an element of $\Sigma_{n+1}$,
so this follows from Lemma \ref{lem:embed} above.

\item Cathelineau's Remark~2.5 points out that, while compatibility
with the $\lambda$-structures is proven for the isomorphism
$K_n(A,I)_{\Q}\cong HC_{n-1}(A,I)$ of \cite[6.2]{OW},
which is the rational homotopy character \eqref{def:relChern},
it agrees with the isomorphism $\rho$ constructed by
Goodwillie in \cite{G86}.
\end{enumerate} 

\section{Appendix: Space-level versions of Cathelineau's Theorem}
\label{app:B} 

Cathelineau's Theorem \ref{thm:cath} concerns the $K$-groups
$K_m(A,I)=\pi_mK(A,I)$ of a nilpotent ideal in a ring.
In this Appendix, we develop a space-level version
(Theorem \ref{thm:KOIcompatible}) of Cathelineau's Theorem. 
That is, for every nilpotent sheaf of ideals $I$,
the rational homotopy character \eqref{eq:KAI-HC} from the simplicial presheaf
$K(\cO,I): U\mapsto K(\cO(U),I(U))$ to $HC(\cO,I)[1]\cong HN(\cO,I)$
preserves the $\lambda$-filtration. By Proposition \ref{ch=rht},
the same is true for the relative Chern character up to homotopy.
Theorem \ref{KiOI}, the spectrum analogue of Theorem \ref{thm:KOIcompatible},
is needed for
the proof of Theorem \ref{Hinf-KI}, which in turn is a key ingredient in
the proof of Theorem \ref{thm:mth}; our Main Theorem \ref{main-intro}
is a special case of \ref{thm:mth}. Theorem \ref{thm:KOIcompatible} is
also used for the scheme-theoretic Theorem \ref{KUIcompatible}.

We will work with a presheaf model \eqref{eq:KOI-HC} of the affine
rational homotopy character \eqref{eq:KAI-HC}
of a sheaf $I$ of nilpotent ideals. 
The remaining steps, which were alluded to
in Step~3 of Appendix \ref{app:A}, follow Cathelineau's
construction in \cite{C91}. Compatible presheaf operations $\Lambda^k$
and $\lambda^k$ are constructed on the $X_n$ and $x_n$ in \eqref{eq:phi}
and \ref{lambda-k-commute}.  The full compatibility is given in
Theorem \ref{thm:KOIcompatible}.

We need to re-introduce some of the notation used in Appendix \ref{app:A}.
For each $n$, any ideal $I$ determines a sheaf of subgroups
$T_n(\cO,I)$ of $GL_n(\cO)$; we form the simplicial subsheaf
$X_n(\cO,I)=\bigcup_{\sigma\in\Sigma_n} BT_n^\sigma(\cO,I)$
of $BGL_n(\cO)$ by conjugating by permutation matrices.
By abuse of notation, we shall write $BGL(\cO)^+$ for the functorial
fibrant model $\Z_\infty BGL(\cO)$ of the plus construction, and write
$K(\cO,I)$ for the homotopy fiber of $BGL(\cO)^+ \to BGL(\cO/I)^+$.
Recall from Theorem \ref{OW6.1} that there is a homology isomorphism
from $X(\cO,I)=\bigcup X_n(\cO,I)$ to $K(\cO,I)$, i.e.,
the natural map $C_*X(\cO,I) \to C_*K(\cO,I)$ is
a quasi-isomorphism when evaluated at any $U$. Therefore
$K_*(\cO,I)\cong\Prim H_*(X(\cO,I))$.

The sheafification of Theorem~\ref{OW2.4} (see \ref{OW2.4bis}) yields a
$\Sigma_n$-equivariant quasi-isomorphism:
$\phi: C_*X_n(\cO,I) \map{\simeq} x_n(\cO,I)$. 
This yields the presheaf version of the rational homotopy character
\eqref{eq:KAI-HC}: 
\begin{equation}\label{eq:KOI-HC}
K(\cO,I) \map{h} C_*K(\cO,I)\ {\buildrel\simeq\over\leftarrow} \
C_*X(\cO,I) \map{\simeq} x(\cO,I) \to HC(\cO,I)[1].
\end{equation}
By Proposition \ref{ch=rht}, the map
$C_*X(\cO,I)\to HC(\cO,I)[1]\simeq HN(\cO,I)$ in \eqref{eq:KOI-HC}
is naturally homotopy to $ch^-$ and hence
\eqref{eq:KOI-HC} is homotopy equivalent to the presheaf form of the
relative Chern character. 
\goodbreak
In order to invert the backwards arrow in \eqref{eq:KOI-HC}, we
use the global projective closed model structure on simplicial presheaves
of sets (discussed in Appendix \ref{app:C});
this is the simplicial closed model structure in which a
presheaf map $f$ is a fibration or a weak equivalence if
$f(U)$ is one for each Zariski open $U$. Let $K(\cO,I)'$ be the
cofibrant replacement for $K(\cO,I)$, and factor the backwards map in
\eqref{eq:KOI-HC} as
$C_*K(\cO,I)\ {\buildrel\simeq\over\twoheadleftarrow}\ C
{\buildrel\simeq\over\leftarrowtail}\ C_*X(\cO,I)$. Then $h$ lifts to a map
$h':K(\cO,I)'\to C$ and, since $x(\cO,I)$ is fibrant,
\eqref{eq:KOI-HC} lifts to a map 
\begin{equation}\label{eq:KOI'-HC}
K(\cO,I)' \map{h'} C \map{\simeq} x(\cO,I) \to HC(\cO,I)[1].
\end{equation}
Since the relative Chern character $ch$ is defined using the same backwards
arrow, 
the remarks after \eqref{eq:KOI-HC} show that
$ch$ lifts to a map $K(\cO,I)'\map{h'}C\to HN(\cO,I)$,
homotopy equivalent to \eqref{eq:KOI'-HC}, followed by
$HC(\cO,I)[1]\simeq HN(\cO,I)$. 

\medskip
Our next step is to construct the operations $\lambda^k$ on
$C_*X_n(\cO,I)$ and $x_n(\cO,I)$.  Following Cathelineau, the exterior
power operations determine group maps $\Lambda_{\times,n}^k$ from
$T_n(\cO,I)$ to $T_{\binom nk}(\cO,I)$ and these induce
maps of simplicial sheaves
\begin{equation}\label{eq:phi}
\Lambda_{\times,n}^k: X_n(\cO,I) \to X_{\binom nk}(\cO,I).
\end{equation}
The exterior power operations also determine Lie algebra maps
$\Lambda_{+,n}^k$ from $\frakt_n(\cO,I)$ to
$\frakt_{\binom nk}(\cO,I)$ and hence chain maps on $x_n(\cO,I)$.
These two constructions are compatible in the following sense.
Let $h:X\to C_*X$ denote the Hurewicz map; 
as in \eqref{eq:KAI-HC}, we regard non-negative
chain complexes as spaces by Dold-Kan.

\begin{lem}\label{lem:Lambda-k}
Let $I$ be a nilpotent sheaf of ideals on $T$, and $n,k$ natural numbers.
The following diagram commutes: 
\begin{equation*}\xymatrix{
X_n(\cO,I) \ar[d]^{\Lambda_{\times,n}^k} \ar[r]^{h} &
C_*X_n(\cO,I) \ar[d]^{\Lambda_{\times,n}^k} \ar[r]_{\phi}^\simeq &
x_n(\cO,I) \ar[d]^{\Lambda_{+,n}^k} \\
X_{\binom nk}(\cO,I) \ar[r]^h &
C_*X_{\binom nk}(\cO,I) \ar[r]_{\phi}^\simeq
& x_{\binom nk}(\cO,I).
}\end{equation*}
\end{lem}

\begin{proof}
This is immediate from the naturality of the  map \eqref{eq:RHT};
see \cite[5.11]{SW}. The left square commutes by naturality of the
Hurewicz map $h$
(see \cite[8.3.9]{WeibelHA94}).
\end{proof}

Now the homology of $C_*X(\cO,I)$ is a Hopf algebra by \ref{OW6.1.1},
and the same is true for $x(\cO,I)$ by \eqref{eq:hocolim}.
Composing $\Lambda^k_n$ with the (chain-level) antipode,
we can define maps $-\Lambda^k_n$ on
$C_*X_n(\cO,I)$ and $x_{n}(\cO,I)$. Combining with the direct sum of
matrices, we can define maps $\lambda_n^k=\oplus \pm\Lambda^i_n$ on
$C_*X_n(\cO,I)$ and $x_{n}(\cO,I)$, taking values in the stable
complexes $C_*X(\cO,I)$ and $x(\cO,I)$. (The number of times the
factor $\pm\Lambda^i_n$ occurs in this expression is suppressed for
legibility.)

In the case of $X(\cO,I)$ we can do better because
$X(\cO,I)^+\map{\sim} K(\cO,I)$ is a weak equivalence by
\ref{OW6.1}. For this we
let $X'_n(\cO,I)$
be a cofibrant replacement for $X_n(\cO,I)$, so that
the maps $\Lambda_{\times,n}^k$ of \eqref{eq:phi} lift to maps
$X'_n(\cO,I)\to K(\cO,I)'$. Using the homotopy inverse on $K(\cO,I)'$
from Proposition \ref{Hinverse}, we also have presheaf maps
$-\Lambda_{\times,n}^k$.  Using the $H$-space structure on $K(\cO,I)'$,
we can combine the maps $\pm\Lambda_{\times,n}^k$ into a 
presheaf map $\lambda^k_{\times,n}: X'_n(\cO,I)\to K(\cO,I)'$ for each $n$.
Given this stabilization, Lemma \ref{lem:Lambda-k}
and \eqref{eq:KOI-HC} imply:

\begin{prop}\label{lambda-k-commute}
The following diagram commutes 
for each $n$ and $k$:
\begin{equation*}\xymatrix{
X'_n(\cO,I) \ar[d]^{\lambda_{\times,n}^k} \ar[rr]^{h} &&
C_*X'_n(\cO,I) \ar[d]^{\lambda_{\times,n}^k} \ar[r]^\simeq &
x_n(\cO,I) \ar[d]^{\lambda_{+,n}^k} \\
K(\cO,I)' \ar[r]^{h} & C_*K(\cO,I)' &  \ar[l]^{\simeq}
C_*X(\cO,I)\ar[r]^\simeq & x(\cO,I).
}\end{equation*}
\end{prop}

This completes the construction of the $\lambda^k$ on $X'_n(\cO,I)$ and
$x_n(\cO,I)$. \smallskip

We may of course pass to the inductive limit over $n$ in the right square
of \eqref{lambda-k-commute}, defining compatible operations on
$C_*X'(\cO,I)$ and $x(\cO,I)$, a passage made possible
because the symmetric group action is homologically trivial by
Lemma \ref{lem:embed}.
(This is observed in \cite[Lemma 3]{C91} and in Step 3(3) of
Appendix \ref{app:A} to this paper.)

We may also perform this stabilization on $X'_n(\cO,I)$. Indeed,
each map $\lambda^k_n$ agrees with the composition
$X'_n(\cO,I)\to X'_{n+1}(\cO,I)\to K(\cO,I)$
up to conjugation by an element of $\Sigma_\infty$. As argued by
Cathelineau on p.~601 of \cite{C91}, using $St(\Z)$ in place of $St(A)$,
there are natural based homotopies between them, so that they are even
base-point homotopic as presheaf maps. Thus they define a map
\begin{equation}\label{eq:XtoK}
X(\cO,I) {\buildrel \sim \over \longleftarrow}
\hocolim X'_n(\cO,I) \map{\lambda^k_n} K(\cO,I).
\end{equation}

The operations $\lambda_{+,n}^k$ in \ref{lambda-k-commute} are
compatible with the operations on cyclic homology.
As observed by Cathelineau on p.~600 of \cite{C91}, the maps
$x_n(\cO,I)\to\wedge^* gl_n(\cO,I)$ and the Loday-Quillen map both commute
with the operations $\lambda_n^k$, so we get commutative diagrams:
\begin{equation}\label{eq:C-p.600}
\xymatrix{
x_n(\cO,I) \ar[d]^{\lambda_{+,n}^k} \ar[r] &
\wedge^*gl_n(\cO,I) \ar[d]^{\lambda_{+,n}^k} \ar[r]
& HC(\cO,I)[1]\ar[d]^{\lambda^k} \\
x(\cO,I)\ar[r] & \wedge^*gl(\cO,I) \ar[r] & HC(\cO,I)[1].
}\end{equation}

The triviality of the symmetric group action in Lemma \ref{lem:embed}
extends to the presheaf level by naturality.
This implies that the maps in \ref{lambda-k-commute} and
\eqref{eq:C-p.600} are compatible (up to chain homotopy) with passage
from $n$ to $n+1$. (This is observed on p.599 of \cite{C91} and in
Step 3(3) of our Appendix \ref{app:A}.)

\medskip\goodbreak
We are now ready to present a space-level
version of Cathelineau's Theorem \ref{thm:cath}:

\begin{thm}\label{thm:KOIcompatible}
The rational homotopy character \eqref{eq:KOI-HC} is compatible with the
operations $\lambda^k$ defined on $K(\cO,I)$ and $HC(\cO,I)$, in the sense
that there is a commutative diagram:

\begin{equation*}
\xymatrix{ K(\cO,I)' \ar[drr]_{\lambda^k} \ar[r]^{\simeq} &
X'(\cO,I)^+ \ar[dr]^{(\ref{eq:XtoK})} 
& \hocolim X'_n(\cO,I)^+ \ar[l]_{\simeq\quad}
\ar[d]^{\lambda_\times^k} \ar[r]^{\quad(\ref{eq:KOI-HC})}
& HC(\cO,I)[1] \ar[d]^{\lambda^k} \\
&& K(\cO,I)' \ar[r]^{(\ref{eq:KOI-HC})} & HC(\cO,I)[1].
}\end{equation*}
The relative Chern character $ch:K(\cO,I)\map{} HN(\cO,I)$
is compatible with the $\lambda$-operations,
up to the natural homotopy of Proposition \ref{ch=rht}.
\end{thm}

\begin{proof}
The diagram \eqref{thm:KOIcompatible} 
commutes because it is obtained by
glueing \eqref{eq:XtoK} together with the homotopy colimits (over $n$)
of the commutative diagrams \eqref{lambda-k-commute} and \eqref{eq:C-p.600},
The upper left map is an inverse to the cofibrant replacement of the
weak equivalence $X(\cO,I)^+\to K(\cO,I)$, constructed using
Lemma \ref{global-we}.

All that remains is to show that
the operations $\lambda^k$ on $K(\cO,I)$ displayed diagonally in
\eqref{thm:KOIcompatible} agree with the usual $\lambda$-operations induced
from the operations on $BGL(\cO)^+$.  To do this, we need to recall
how the usual operations are defined.

As observed in \cite[2.1]{C91}, 
the construction of the usual operations $\lambda^k$ 
may be chosen to follow the pattern described above. Briefly,
one defines natural group maps $\Lambda^k_n$ on $GL_n(\cO)$ as
in \eqref{eq:phi} and applies \ref{Hinverse} to the $H$-space
structure on $BGL(\cO)^+$ to obtain an $H$-inverse
$\iota$ on a cofibrant replacement of $BGL(\cO)^+$,
and then construct maps $\lambda^k_n$ on
$BGL'_n(\cO)$, as described before \ref{lambda-k-commute}.
As in \eqref{eq:XtoK}, this induces a map $\lambda^k$ on $BGL'(\cO)$ and
$BGL(\cO)^+$; finally this induces the relative operations $\lambda^k$
on $K(\cO,I)$. (Compare with \cite[p.511]{Soule}.)
By construction, the inclusions $X'_n(\cO,I)\subset BGL'_n(\cO)$
are compatible with the operations $\lambda^k_n$ just defined, and the
$\lambda^k$ defined in \eqref{eq:XtoK}. Thus we have a commutative diagram:
\addtocounter{equation}{-1}
\begin{subequations}
\begin{equation}\label{eq:C-p.602}
\xymatrix{ {\hskip-12pt}
X'(\cO,I) \ar[dr]_{(\ref{eq:XtoK})}&
\hocolim X'_n(\cO,I) \ar[l]_{\simeq\qquad} \ar[d]^{\lambda_\times^k} \ar[r] &
\hocolim BGL'_n(\cO) \ar[d]^{\lambda_+^k}
\rlap{$\map{\simeq} BGL(\cO)$} \\ 
& K(\cO,I)' \ar[r] & BGL(\cO)^+. 
}\end{equation}
\end{subequations}

\goodbreak
Applying the functorial $+$-construction $\Z_\infty$ converts $X'(\cO,I)$
and $BGL(\cO)$ into $K(\cO,I)$ and $BGL(\cO)^+$ in \eqref{eq:C-p.602},
up to weak equivalence.
As argued by Cathelineau on p.~602 of \cite{C91}, the horizontal composites
$K(\cO,I)\to BGL(\cO)^+$ are part of a presheaf of fibration sequences,
with third term $BGL(\cO/I)^+$. It follows
that the map \eqref{eq:XtoK} does indeed induce the
usual operation $\lambda^k$ on $K(\cO,I)$, as claimed.
\end{proof}

\begin{subremark}\label{psi-compatible}
Since the relative Chern character is multiplicative,
it follows from \ref{thm:KOIcompatible} 
that there is a similar compatibility result (up to homotopy)
for polynomials in the $\lambda^i$.
For example, the Adams operations $\psi^k$ are defined as polynomials;
by \cite[5.3]{Kratzer}, $\psi^k = (-1)^{k-1}k\lambda^k$ on
$\pi_* K(\cO,I)$ and $HC_*(\cO,I)$. Hence the Adams operations
are compatible with the relative Chern character up to natural homotopy.

\end{subremark}
\hskip-12pt\goodbreak

In order to extend Theorem \ref{thm:KOIcompatible} to presheaves of spectra,
recall that we can regard chain complexes as spectra, and morphisms of chain
complexes as morphisms of spectra (see \cite[p.\ 552]{WeibelNil}).
In this way, $HC(\cO,I)[1]$ is a presheaf of (connective) spectra
and the $\lambda^k$ are spectrum endomorphisms of $HC(\cO,I)[1]$.

In a parallel abuse of notation, we shall write $K(\cO,I)$ for the
presheaf of spectra whose initial space is the space $K(\cO,I)$.
Note that the {\em a priori} non-connective $K$-theory spectrum obtained
by evaluation at any $U$ just happens to be connective because it
is the relative $K$-theory spectrum of a nilpotent ideal.
It is easy to see that the rational homotopy character and the
relative Chern character $ch:K(\cO,I)\map{} HN(\cO,I)$ are both
morphisms of spectra, and that $ch$ is a homotopy equivalence
of spectra by Goodwillie's Theorem \cite{G86}.  (See \cite[Exercise
9.10]{TT}; the key is that $ch$ is multiplicative by \cite[\S5]{HJ}).
Using the main result of \cite{cw-agree}, it is easy to see that these
morphisms are homotopic.
Combining these facts with Theorem \ref{thm:KOIcompatible} and Remark
\ref{psi-compatible}, we see that we have proven:
\hskip-5pt

\begin{cor}\label{deloop-lambda}
When $I$ is a nilpotent sheaf of ideals, 
$\lambda^k$ and $\psi^k$ are morphisms of presheaves of spectra,
from  $K(\cO,I)'$ to $K(\cO,I)$, and commute with $ch$ up to
natural homotopy equivalence.
\end{cor}

We conclude with an interpretation in terms of the eigen-components
with respect to these operations. We fix $k\ge2$, and define
$K^{(i)}(\cO,I)$ to be the cofibrant homotopy fiber of
$\lambda^k+(-1)^kk^{i-1}: K(\cO,I)'\to K(\cO,I)$. We leave it as an exercise
to see that (up to homotopy equivalence) this is independent of the
choice of $k>1$, and that  our choice here of $\lambda^k$ rather than
$\psi^k$ is immaterial (use Remark \ref{psi-compatible}).

Recall too from \cite[4.6.7 and 4.5.16]{LodayHC92} that the chain
complex for cyclic homology breaks up into the direct product of
subcomplexes $HC^{(i)}$ on which $\lambda^k+(-1)^kk^i$ and $\psi^k-k^{i+1}$
are acyclic. Passing
to the associated Eilenberg-Mac\,Lane spectra, this means that
$HC^{(i)}(\cO,I)$ is the homotopy fiber of
$\lambda^k+(-1)^kk^i$ 
on $HC(\cO,I)$. A similar description holds for $HN(\cO,I)$; see
\cite[5.1.20]{LodayHC92}. Thus we have proven:
\hskip-5pt

\begin{thm}\label{KiOI}
The rational homotopy character \eqref{eq:KOI-HC} induces maps on each
eigen-component $K^{(i)}(\cO,I)$, fitting into a commutative square:

$$\begin{CD}
K^{(i)}(\cO,I) @>{\simeq}>> HC^{(i-1)}(\cO,I)[1]
	@>{\simeq}>> HN^{(i)}(\cO,I) \\
@VVV   @VVV @VVV \\
K(\cO,I) @>{\simeq}>> HC(\cO,I)[1] @>{\simeq}>> HN(\cO,I).
\end{CD}$$
\end{thm}
\hskip-5pt
Using the fact that $HC(\cO,I)\cong\prod_{i\ge0} HC^{(i)}(\cO,I)$,
we see from \ref{KiOI} that the maps $K^{(i)}(\cO,I) \to K(\cO,I)'$
are split, and that we have
a homotopy equivalence $K(\cO,I)'\map{\simeq}\prod_{i\ge1} K^{(i)}(\cO,I)$.
Thus we deduce:
\hskip-5pt

\begin{cor}\label{prodKiOI}
There is a homotopy commutative diagram:
\begin{equation*}
\minCDarrowwidth10pt
\begin{CD}
K(\cO,I) @>{\simeq}>{ch}> HC(\cO,I)[1] @>{\simeq}>> HN(\cO,I) \\
 @VV{\simeq}V @VV{\simeq}V @VV{\simeq}V \\
\prod_{i=1}^\infty K^{(i)}(\cO,I)  @>{\simeq}>{ch}>
 \prod_{i=1}^\infty HC^{(i-1)}(\cO,I)[1] @>{\simeq}>>
\prod_{i=1}^\infty HN^{(i)}(\cO,I).
\end{CD}\end{equation*}
\end{cor}

\bigskip\goodbreak
\section{Appendix: Simplicial presheaves of Sets}\label{app:C}

In this Appendix we prove some elementary results about simplicial presheaves
of sets, 
which are used in Appendix \ref{app:B}. By a simplicial presheaf we mean a
contravariant functor from a fixed small category $\mathbb T$
to the category of simplicial sets.

\smallskip

\begin{lem}\label{global-we}
Let $f:E\to E'$ be a weak equivalence 
in any model category, with $E$ fibrant and $E'$ cofibrant.
Then there is a map $g:E'\to E$ such that $gf$ and $fg$ are
weak equivalences.
\end{lem}

\begin{proof}
Consider the mapping cylinder $T$ of $f$, weak equivalent to $E'$,
constructed using the (CM5) factorization of $E\vee E'\to E'$ in the
model structure.  
Since $E'$ is cofibrant, $i:E\to T$ is a trivial cofibration.
Since $E$ is fibrant, $i$ has a retract by axiom (CM4).
But then the map $g:E'\to T\to E$ has the desired properties;
the homotopy from the identity on $E$ (which is $E\to T\to E$) to $gf$
is given by the map $T\map{\sim}E'\map{\sim}T$; the homotopy from $fg$ to
the identity on $E'$ is similar.
\end{proof}

We say that a map $E\to E'$ of simplicial presheaves is a (global)
weak equivalence if 
each $E(U)\to E'(U)$ is a homotopy equivalence of simplicial sets
(or, more formally, a homotopy equivalence between the corresponding
topological spaces).
A {\it global model structure} on simplicial presheaves is one whose
weak equivalences are the global weak equivalences.

\begin{prop}\label{Hinverse}
Let $H$ be a fibrant presheaf of connected $H$-spaces, and
$H'\buildrel{\sim}\over\twoheadrightarrow H$ its cofibrant replacement
in a global model structure.
Then there is a presheaf map $\iota:H'\to H'$ such that each $H'(U)$ is an
$H$-group with homotopy inverse $\iota(U)$.
\end{prop}

\begin{proof}
Consider the shear map $\phi:H\times H\to H\times H$ defined by
$\phi(x,y)=(x,xy)$. By \cite[X.2.1]{Wh}, each $\phi(U)$ is a
homotopy equivalence, i.e., $\phi$ is a global weak equivalence.
If $C\buildrel{\sim}\over\twoheadrightarrow H\times H$
is a cofibrant replacement, then $C$ is fibrant (because $H\times H$ is)
and the shear map lifts to a map $\phi: C\to C$.
By \ref{global-we}, there is a map $\psi:C\to C$ which is
a homotopy inverse of $\phi$ on each $U$. Let
$i_1:H'\to C$ and $pr_2:C\to H'$ be the cofibrant replacements of the
inclusion $H=H\times\ast \to H\times H$ and the projection $H\times
H\to \ast\times H=H$.  By \cite[III.4.17]{Wh}, each $H'(U)$ is an
$H$-group, via $H'(U)\map{\sim}H(U)$; the composite
$H' \map{i_1} C \map{\psi} C \map{pr_2} H'$
will be the desired map $\iota$ which is a homotopy inverse for each
$H$-space $H'(U)$.
\end{proof}

We can do better using the {\it global projective} closed model structure,
in which fibrations are defined objectwise; this simplicial closed model
structure was introduced by Quillen in \cite[Ch.II, \S4, Th.4]{QHA}.
Indeed, the representable presheaves $h_X:U\mapsto\Hom_{\mathbb T}(U,X)$
form a family of small projective generators; see \cite[pp. 314]{BK}.
It follows that, for any $X$ in $\mathbb T$ and any cofibration
of simplicial sets $K\rightarrowtail L$,
the canonical arrow $h_X\otimes K\to h_X\otimes L$ is a cofibration,

The following result does not seem to be in the literature.

\begin{thm}\label{prod-cofibrant}
Assume that $\mathbb T$ has finite products.
For the global projective closed model structure,
the product of cofibrant presheaves is cofibrant.
\end{thm}

Recall from the proof of \cite[Ch.II, \S4, Th.4]{QHA} that
every simplicial presheaf has a canonical cofibrant replacement
$C^{(\infty)}$, in which the
``$n$-skeleton'' $C^{(n)}$ is obtained by attaching cells to $C^{(n-1)}$
of the form $h_X\otimes\Delta^n$ along
maps $h_X\otimes\partial\Delta^n\to C^{(n-1)}$.

\begin{proof}
It suffices to assume that $C_1$, $C_2$ are canonical cofibrant presheaves
and show that $C_1\times C_2$ is cofibrant. Since
$h_X\times h_Y=h_{X\times Y}$, the result is true for representable
presheaves, and in particular for the $0$-skeleton of $C_1\times C_2$.
Inductively, suppose that $C^n=C_1^{(n)}\times C_2^{(n)}$ is cofibrant
and that $C_i^{(n+1)}$ is obtained from $C_i^{(n)}$ by pushing out along
coproducts of cells of the form $A_i=
h_{X_i}\otimes\partial\Delta^n\rightarrowtail h_{X_i}\otimes\Delta^n
=B_i$. Then
$$ 
C^n \to C^{(n+1)}\times C_2^{(n)} \to
C^{(n+1)}\times C_2^{(n)} \cup_{C^n} 
C^{(n)}\times C_2^{(n+1)} \to C^{(n+1)}\times C_2^{(n+1)}
$$
is a cofibration since the first two maps are just pushouts along
co-base extensions of
$A_i\to B_i$ ($i=1,2$), and the final map is a co-base extension of
$A_1\times A_2\rightarrowtail B_1\times B_2$,
which is a cofibration because it is a coproduct of maps
$h_{X_i\times X_j}\otimes (K\rightarrowtail L)$.
\end{proof}

\begin{cor}
If $H$ is a presheaf of connected fibrant $H$-spaces, its cofibrant
replacement $H'$ (in the global projective closed model structure) is
a presheaf of $H$-groups.
\end{cor}

\begin{proof}
Since $H'\times H'$ is cofibrant by \ref{prod-cofibrant},
the product $H\times H\to H$ lifts to a product $H'\times H'\to H'$,
and the homotopy inverse $\iota:H'\to H'$ exists by \ref{Hinverse}.
\end{proof}

\goodbreak
\subsection*{Acknowledgements}
The authors would like to thank Rick Jardine for his help with
fibrant replacements, especially in the proof of Theorem \ref{thm:fibrant}.
We would also like to thank the referee for pointing out that the
relative Chern character and the rational homotopy characters needed to
be identified, a problem solved in \cite{cw-agree}.

\end{document}